\documentclass[11pt]{article}
\usepackage{amsmath,amssymb}
\usepackage{graphicx} 

\newcommand{\cqfd}{\hfill $\square$}
\newcommand{\mino}[2]{\min(#1, #2)}

\newcommand{\Nn}{\mathbb{N}}

\newcommand{\pla}{\ensuremath{A}}
\newcommand{\plb}{\ensuremath{B}}

\newcommand{\proof}{\noindent{\bf Proof:} }

\usepackage{version}

\newtheorem{lemma}{Lemma}[section]
\newtheorem{proposition}{Proposition}[section]
\newtheorem{theorem}{Theorem}[section]
\newtheorem{corollary}{Corollary}[section]
\newtheorem{remark}{Remark}[section]
\newtheorem{definition}{Definition}[section]

\begin{document}

\title{A Stochastic Analysis of some Two-Person Sports}

\author{Davy \sc{Paindaveine} \footnote{E.C.A.R.E.S. and  D\'epartement de Math\'ematique, av. F.D. Roosevelt 50 - CP114, B-1050 Brussels}
              Yvik \sc{Swan}\footnote{D\'epartement de Math\'ematique,  Campus Plaine - CP 210, B-1050 Brussels}}
\maketitle
\begin{abstract}
We consider two-person sports where each rally is initiated by a \emph{server}, the other player (the \emph{receiver}) becoming the server when he/she wins a rally.  Historically, these sports  used a  scoring based on the \emph{side-out scoring system},  in which points are only scored by the server. Recently, however, some federations have switched to the \emph{rally-point scoring system} in which a point is scored on every rally. As various authors before us, we study how much this change affects the game. Our approach is based on a \emph{rally-level analysis} of the process through which, besides the  well-known probability distribution of the scores, we also obtain the distribution of the number of rallies.  This  yields a comprehensive knowledge of the process at hand, and  allows for an in-depth comparison of both scoring systems. In particular, our results {help} to explain why the transition from one scoring system to the other  has more important implications than those predicted from game-winning probabilities alone. Some of our findings are quite surprising, and unattainable through Monte Carlo experiments. Our results are of high practical relevance to  international federations and local tournament organizers alike, and also open the way to efficient estimation of the rally-winning probabilities, which should have a significant impact on the quality of ranking procedures.
\end{abstract}

\noindent {\small \textbf{Keywords.}  
{Combinatorial derivations;
Duration analysis;
Point estimation;
Ranking procedures;
Scoring rules;
Two-person sports
} 



\section{Introduction.} \label{intro}
We consider a class of two-person
sports for which each rally is initiated by a \emph{server}---the other player is then called the \emph{receiver}---and for which the rules and scoring system satisfy one of the following two definitions. 
\begin{quote} \emph{Side-out scoring system}: (i) the server in the first rally is determined by flipping a coin. (ii) If a rally is won by the server, the latter scores a point and serves in the next rally. Otherwise, the receiver becomes the server in the next rally, but no point is scored. (iii) The winner of  the \emph{game} is the first player to score $n$ points. 
\end{quote}
\begin{quote}   \emph{Rally-point scoring  system}: (i) the server in the first rally is determined by flipping a coin. (ii) If a rally is won by the server, the latter  serves in the next rally. Otherwise, the receiver becomes the server in the next rally. A point  is scored after each rally. (iii) The winner of  the \emph{game} is the first player to score $n$ points. 

\end{quote}

A \emph{match} would typically consist of a sequence of such games, and the winner of the match is the first player to win~$M$ games. Actually, it is usually so that  in game~$m\geq 2$, the first server is not determined by flipping a coin, but rather according to some prespecified rule: the most common one states that the first server in game~$m$ is the winner in game~$m-1$, but alternatively, the players might simply take turns as the first server in each game until the match is over. It turns out that, in the probabilistic model we consider below, the probability that a fixed player wins the match is the same under both rules; see {Anderson (1977).} This clearly allows us to focus on a single game in the sequel---as in most previous works in the field (references will be given below). Extensions of our results to the match level can then trivially be obtained by appropriate conditioning arguments, taking into account the very rule adopted for determining the first server in each game.

The side-out scoring system has been used in various sports, sometimes up to tiny unimportant refinements, {involving typically, in case of a tie at $n-1$, the possibility (for the receiver)  to choose whether the game should be played to $n+\ell$ (for some fixed~$\ell\geq 2$) or to~$n$}; see Section~\ref{probagame}. When based on the so-called English scoring system, Squash currently uses $(n,M)=(9,3)$. 
Racquetball is essentially characterized by~$(n,M)=(15,2)$ (the possible third game is actually played to 11 only). 
Until 2006, Badminton was using $(n,M)=(15,2)$ and $(n,M)=(11,2)$ for men's and women's singles, respectively---with an exception in 2002, where $(n,M)=(7,3)$ was experienced. 
Volleyball, for which the term \emph{persons} above should of course be understood as \emph{teams}, 
was based on $(n,M)=(15,3)$ until 2000. 
In both badminton and volleyball, this scoring system was then replaced with the {rally-point system}. Similarly, squash, at the international level, now is based on the American version of its scoring system, which is nothing but the rally-point system, in this case with~$(n,M)=(11,3)$. Investigating the deep implications of this transition from the original side-out scoring system to the rally-point scoring system was one of the main motivations  for this work; see Section~\ref{amerisec}. 

Irrespective of the scoring system adopted, the most common probabilistic model for the sequence of rallies  assumes that the rally outcomes are mutually independent and are, conditionally {on} the server, identically distributed. This implies that the game is governed by the parameter~$(p_a,p_b)\in[0,1]\times[0,1]$, where the \emph{rally-winning probability} $p_a$ (resp., $p_b$) is the probability that Player $\pla$ (resp., Player $\plb$) wins a rally when serving. This means that players do not get tired through the match, or that they do not get nervous when playing crucial points. More precisely, they might get tired or nervous, but if they do, they should do so at the same moment and to the same extent, so that it does not affect their respective rally-winning probability. We will throughout refer to this probabilistic model as the \emph{server model}, in contrast with the \emph{no-server model} in which any rally is won by~$\pla$ with probability~$p$ irrespective of the server, that is, the submodel obtained when taking $p=p_a=1-p_b$.

The probabilistic properties of a single game {played under the side-out scoring system have} been  investigated in various works. Hsi and Burich~(1971) attempted to derive the probability distribution of game scores---in the sequel, we simply speak of the \emph{score distribution}---in terms of~$p_a$ and~$p_b$, but their derivation based on standard combinatorial arguments {was} wrong. The correct score distribution (hence also the resulting game-winning probabilities) was first obtained in Phillips~(1978) by applying results on sums of random variables having the modified geometric distribution. Keller~(1984)  derived probabilities of very extreme scores, whereas Marcus~(1985) derived the complete score distribution in the no-server model. Strauss and Arnold (1987), by identifying the point earning process as a Markov chain, obtained more directly the same general result as  Phillips~(1978). They further used the score distribution to define maximum likelihood estimators and moment estimators of the rally-winning probabilities (both in the server and no-server models), and based on these estimates a ranking system (relying on Bradley-Terry paired comparison methods) for the players of a league or tournament. 
Simmons~(1989) determined the score distribution under the two scoring systems, this time by using a quick and direct combinatorial analysis of a single game. He discussed handicapping and strategies (for deciding whether {the receiver should go for a game played to~$n+\ell$ or not in case of a tie at~$n-1$}), and attempted a comparison of the two scoring systems. More recently, Percy (2009) used Monte Carlo simulations to compare game-winning probabilities and \emph{expected} durations for both scoring systems in the no-server model. 

To sum up, the score distributions have been obtained through several different probabilistic methods, and were used to discuss several aspects of the game. In contrast, the distribution of the number of rallies  needed to complete a single game~($D$, say) remains virtually unexplored for the side-out scoring system (for the rally-point scoring system, the distribution of~$D$ is simply determined by the score distribution). To the best of our knowledge, the only theoretical result on~$D$ under the side-out scoring system provides lower and upper bounds for the \emph{expected value} of~$D$; see~(20) in Simmons (1989), or~(\ref{simmons}) below. Beyond the lack of exact results on~$D$ (only approximate theoretical results or simulation-based results are available so far), it should be noted that only the expected value of~$D$ has been studied in the literature. This is all the more surprising because, in various sports (e.g., in badminton and volleyball), uncertainty about~$D$---which is related to its \emph{variance}, not to its expected value---was one of the most important arguments to switch from the side-out scoring system to the rally-point scoring system. Exact results on the moments of~$D$---or even better, its distribution---are then much desirable as they would allow to investigate whether the transition to the rally-point system indeed reduced uncertainty about~$D$. More generally, precise results on the distribution of~$D$ would allow for a much deeper comparison of both scoring systems. They would also be of high practical relevance, e.g.,  to tournament organizers, who need planning their events and deciding in advance the number of matches---hence the number of players---the events will be able to host. 

For the side-out scoring system, however, results on the distribution of~$D$ cannot be obtained from a \emph{point-level} analysis of the game. That is the reason why the present work rather relies on a  \emph{rally-level}  combinatorial analysis. This allows to get of rid of the uncertainty about the number of rallies  needed to score a single point, and results into an exact computation of the distribution of~$D$---and actually, even of the  number of rallies  needed to achieve any fixed score. We derive explicitly the expectation and variance of~$D$, and use our results to compare the two scoring systems  not only in terms of game-winning probabilities, but also in terms of durations. Our results reveal significant differences between both scoring systems, and help to explain why the transition from one scoring system to the other has more important implications than those predicted from game-winning probabilities alone. 
As suggested above, they could be used by tournament organizers to plan accurately their events, but also by national or international federations to better perform the possible transition from the side-out scoring system to the rally-point one; see Section~\ref{finalcom} for a discussion. Finally, our results open the way to efficient estimation of the rally-winning probabilities (based on observed scores and durations), which might have important consequences for the resulting ranking procedures, since rankings usually are to be based on small numbers of ``observations" (here, games).

The outline of the paper is as follows. In Section~\ref{probagame}, we describe our rally-level analysis of a single game played under the side-out scoring system, and show that it also leads to the score distribution already derived in Phillips~(1978), Strauss and Arnold (1987), and Simmons (1989). Section~\ref{duration} explains how this rally-level analysis further provides (i) the expectation and variance of the number of rallies  needed to achieve a fixed score (Section~\ref{expectgame}) and also (ii) the corresponding exact distribution (Section~\ref{distriduration}). In Section~\ref{amerisec}, we then use our results in order to compare the side-out and rally-point scoring systems, both in terms of game-winning probabilities (Section~\ref{probacompare}) and durations (Section~\ref{durationcompare}). In Section~\ref{simu}, we perform  Monte Carlo simulations and compare the results with our theoretical findings. Section~\ref{finalcom} presents the conclusion and provides some final comments. Finally, an appendix collects proofs of technical results. 
\section{Rally-level derivation of the score distribution under the side-out scoring system.} \label{probagame}

In this section, we conduct our rally-level analysis of a single game played under the side-out scoring system. We will make the distinction between $\pla$-games and $\plb$-games, with the former (resp., the latter) being defined as games in which Player~$\pla$ (resp., Player $\plb$) is the first server. Wherever possible, we will state our results/definitions in the context of $\pla$-games only; in such cases  the corresponding  results/definitions for $\plb$-games can then be obtained by exchanging the roles played by~$\pla$ and~$\plb$, that is, by exchanging (i)~$p_a$ and $p_b$ and~(ii) the number of points scored by each player. Whenever not specified, the server~$S$ will be considered random, and we will denote by~$s_a:={\rm P}[S=A]$ and~$s_b:={\rm P}[S=B]=1-s_a$ the probabilities that the game considered is an $A$-game and a $B$-game, respectively. This both covers  games where the first server is determined by flipping a coin and games where the first server is fixed (by letting~$s_a\in\{0,1\}$).

Our rally-level analysis of the game will be based on the concepts of \emph{interruptions} and \emph{exchanges} first introduced  in Hsi and Burich~(1971). More precisely, we adopt

\begin{definition} \label{defone}
An \emph{$\pla$-interruption} is a sequence of rallies  in which $\plb$ gains the right to serve from $\pla$, scores at least one point, then (unless the game is over) relinquishes the service back to $\pla$, who will score at least one point. An \emph{exchange} is a sequence of two rallies  in which one player gains the right to serve, but immediately loses this right before he/she scores any point (so that the potential of consecutive scoring by his/her opponent is not interrupted).
\end{definition}
\vspace{1mm}

We point out that $\pla$-interruptions are characterized in terms of score changes only (and in particular may contain one or several exchanges) and that, at any time, an exchange clearly occurs with probability~$q:=q_a q_b:=(1-p_a)(1-p_b)$. 






Now, for~$C\in \{A,B\}$, denote by 
$
E^{\alpha,\beta,C}(r,j)
$
the event associated with a sequence of rallies  that (i) gives raise to $\alpha$ points  scored by Player~A and $\beta$ points scored by Player~B, (ii)  involves exactly~$r$ $\pla$-interruptions and $j$ exchanges, and (iii) is such that Player~$C$ scores a point in the last rally; the superscript~$C$ therefore indicates who is scoring the last point, and it is  assumed here that $\alpha>0$ (resp., $\beta>0$) if $C=A$ (resp., if $C=B$). We will write 
$$
p_{C_1}^{\alpha,\beta,C_2}(r,j) := {\rm P}[ E^{\alpha,\beta,C_2}(r,j) | S=C_1 ],
\qquad
C_1,C_2\in \{A,B\}.
$$
We then have the following result
 (see the Appendix for the proof).

\begin{lemma} \label{probascoreinterexch}
Let $\gamma_0 := \min\{\beta,1\}$, ${\gamma_1} := \min\{\alpha, \beta\}$, and $\gamma_2  := \min\{\alpha, \beta-1\}$. Then, setting $\binom{-1}{-1}:=1$, we have $
\textstyle
p_{\pla}^{\alpha,\beta,\pla}(r,j)
=
 {\alpha+\beta+j-1 \choose j}  {\alpha \choose r} {\beta-1 \choose r-1} p_a^{\alpha} p_b^{\beta} q^{r+j}$,
 $r\in\{\gamma_0,\ldots, \gamma_1\}$,
$j\in\Nn,
$
and 
$
\textstyle
p_{\pla}^{\alpha,\beta,\plb}(r,j)
=
 {\alpha+\beta+j-1 \choose j} {\alpha \choose r-1} {\beta-1 \choose r-1} p_a^{\alpha} p_b^{\beta} q_a q^{r+j-1},$
$r\in\{ 1,\ldots, \gamma_2+1 \},$
$j\in\Nn$.
\end{lemma}
\vspace{2mm}


By taking into account all possible values for the numbers of $\pla$-interruptions and exchanges, Lemma~\ref{probascoreinterexch} quite easily leads to the following result (see the Appendix for the proof), which then trivially provides the score distribution in an $\pla$-game, hence also the corresponding game-winning probabilities.



\vspace{1mm}
\begin{theorem} \label{probascore}
Let $p_{C_1}^{\alpha,\beta,C_2}:={\rm P}[E^{\alpha,\beta,C_2}|S=C_1]$, where~$E^{\alpha,\beta,C_2}:=\cup_{r,j} \, E^{\alpha,\beta,C_2}(r,j)$, with~$C_1,C_2\in \{A,B\}$. Then
$
\textstyle
p_{\pla}^{\alpha,\beta,\pla}
=
\frac{p_a^{\alpha} p_b^{\beta}}{(1-q)^{\alpha+\beta}}
\sum_{r=\gamma_0}^{\gamma_1} 
{\alpha \choose r} {\beta-1 \choose r-1}  q^{r}
$
and
$\textstyle
p_{\pla}^{\alpha,\beta,\plb}= $\linebreak
$\frac{p_a^{\alpha} p_b^{\beta}q_a}{(1-q)^{\alpha+\beta}}
 \sum_{r=1}^{\gamma_2+1} 
{\alpha \choose r-1} {\beta-1 \choose r-1}  q^{r-1} 
$.
\end{theorem}
\vspace{1mm}

In the sequel, we denote game scores by couples of integers, where the first entry (resp., second entry) stands for the number of points scored by Player~$A$ (resp., by Player~$B$). 
With this notation, a $C$-game ends on the score~$(n,k)$ (resp., $(k,n)$), $k\in\{0,1,\ldots,n-1\}$, with probability~$p_{C}^{n,k,\pla}$  (resp.,~$p_{C}^{k,n,\plb}$), hence is won by~$A$ (resp., by~$B$) with the (\emph{game-winning}) probability
$$
p_C^{\pla}
:={\rm P}[E^A | S=C]
=\sum_{k=0}^{n-1} p_C^{n,k,\pla}$$ 
(resp., $p_C^{\plb}:=1-p_C^{\pla}$); throughout, $E^A := \cup_{k=0}^{n-1} E^{n,k,A}$ (resp., $E^B := \cup_{k=0}^{n-1} E^{k,n,B}$) denotes the event that the game---irrespective of the initial server---is won by~$A$ (resp., by~$B$). Of course, 
unconditional on the initial server, we have
$$
p^{n,k,A}:={\rm P}[E^{n,k,A}]=p_A^{n,k,A} s_a+ p_B^{n,k,A} s_b,
\qquad
p^{k,n,B}:={\rm P}[E^{k,n,B}]=p_A^{k,n,B} s_a+ p_B^{k,n,B} s_b,
$$
and
$$
p^{C}:={\rm P}[E^C]=p_A^Cs_a+p_B^Cs_b,
$$ 
for $C\in\{A,B\}$.

Figures~\ref{fig1}(a)-(b)  present, for $n=15$, the score distributions associated with~$(p_a,p_b)=(.7,.5)$, $(.6,.5)$, $(.5,.5)$, and $(.4,.5)$. We reversed the $k$-axis in Figure~\ref{fig1}(b), since, among all scores associated with a victory of~$\plb$, the score (14,15) can be considered the closest to the score (15,14) (associated with a victory of~$A$). It then makes sense to regard Figures~\ref{fig1}(a)-(b) as a single plot. The resulting ``global" probability curves are {quite smooth and, as expected, unimodal (with the exception of the $p_a=p_b=.5$ curve, which is slightly bimodal)}. It appears that these score distributions are extremely sensitive to $(p_a,p_b)$, as are the corresponding game-winning probabilities ($p_\pla^\pla$ ranges from .94 to .22, when, for fixed~$p_b=.5$, $p_a$ goes from $.7$ to~$.4$).
For $p_a=p_b=.5$, we would expect the global probability curve to be symmetric. The advantage Player~$A$ is given by serving {first} in the game, however, makes this curve slightly asymmetric; this is quantified by the corresponding probability that $\pla$ wins the game, namely~$p_\pla^\pla=.53>.47=p_\pla^\plb$. 

As mentioned in the Introduction, sports based on the side-out scoring system may involve tie-breaks in case of a tie at~$n-1$. This means that, at this tie, the receiver
has the option of playing through to~$n$ or  ``setting to~$\ell$" (for a  fixed $\ell\geq 2$), in which case the winner is the first player to score~$\ell$ further points. For instance,  games in the current side-out scoring system for squash are played to~$n=9$ points, and the receiver, at~$(8,8)$, may decide whether the game is to~$9$ or $10$ points ($\ell=2$). Before the transition to the rally-point system in 2006, similar tie-breaks were used in badminton, there with~$n=15$ and~$\ell=3$. Assuming that the game is always set to~$\ell$ in case of tie at~$n-1$,
the resulting score distribution can then be easily derived from Theorem~\ref{probascore} by appropriate conditioning; for instance, the score~$(n+\ell-1,n+k-1)$, $k\in\{0,1,\ldots,\ell-1\}$ occurs in an $A$-game with probability~$p_{\pla}^{n-1,n-1,\pla}p_{\pla}^{\ell,k,\pla}+p_{\pla}^{n-1,n-1,\plb}p_{\plb}^{\ell,k,\pla}$. We stress that all results we derive in the later sections can also be extended to scoring systems involving tie-breaks, again by appropriate conditioning. Finally, various papers discuss tie-break strategies (whether to  play through or to set the game to~$\ell$) on the basis of~$p_a$ and~$p_b$; see, e.g., Renick (1976, 1977), Simmons~(1989), or Percy~(2009).



\section{Distribution of the number of rallies  under the side-out scoring system.}
\label{duration}
\vspace{2mm}

As mentioned in the Introduction, the literature contains few results about the number of rallies  $D$ to complete a single game played under the side-out scoring system. Of course, the distribution of~$D$ can always be investigated by simulations; see, e.g., Percy~(2009), where Monte Carlo methods are used to estimate the expectation of~$D$ for a broad range of rally-winning probabilities in the no-server model. To the best of our knowledge, the only available \emph{theoretical} result is due to Simmons~(1989), and provides lower and upper bounds on the \emph{expectation} of~$D$ in an $A$-game conditional on a victory of~$A$ on the score~$(n,k)$. More specifically, letting
\begin{equation} \label{defexp}
e_{C_1}^{\alpha,\beta,C_2}:={\rm E}[ D\, |\, E^{\alpha,\beta,C_2}, S=C_1 ],
\qquad C_1,C_2\in\{ A,B \},
\end{equation}
Simmons' result states that 
\begin{equation}\label{simmons}
\textstyle
(n+k)\,\frac{1+q}{1-q} 
\leq 
e_{A}^{n,k,A}
\leq 
(n+k)\,\frac{1+q}{1-q} + 2 k,\quad k=0,1,\ldots,n-1.
\end{equation}
Unless a shutout is considered (that is, $k=0$), this is only an approximate result, whose accuracy quickly decreases with~$k$. Again, the reason why no exact results are available is that all analyses of the game in the literature are of a \emph{point-level} nature. In sharp contrast, our \emph{rally-level} analysis allows, {\sl inter alia}, for obtaining exact values of all moments of~$D$, as well as its complete distribution.

\subsection{Moments.}
\label{expectgame}
\vspace{2mm}

We first introduce the following notation. Let~$R_{\pla}^{\alpha,\beta,A}$ (resp.,~$R_{\pla}^{\alpha,\beta,B}$) be a random variable assuming {values}~$r=\gamma_0,\gamma_0+1,\ldots,\gamma_1$ 
(resp., $r=1,2,\ldots,\gamma_2 +1$) with corresponding {probabilities}~$
W^{\alpha,\beta,\pla}_\pla(q,r) 
:=
\binom{\alpha}{r}\binom{\beta-1}{r-1}q^r
/\,[\sum_{s=\gamma_0}^{\gamma_1}\binom{\alpha}{s}\binom{\beta-1}{s-1}q^s]
 $ (resp.,  \linebreak$
 W^{\alpha,\beta,\plb}_\pla(q,r) := 
 \binom{\alpha}{r-1}\binom{\beta-1}{r-1}q^{r-1}/\,[\sum_{s=1}^{\gamma_2 +1}\binom{\alpha}{s-1}\binom{\beta-1}{s-1}q^{s-1}]
$). 
Conditioning with respect to the number of $\pla$-interruptions and exchanges then yields the following result (see the Appendix for the proof).

\begin{theorem} \label{momentgenerating}
Let $t\mapsto M_{C_1}^{\alpha,\beta,C_2}(t)={\rm E}[ e^{tD}\, |\, E^{\alpha,\beta,C_2}, S=C_1 ]$, $C_1,C_2\in\{ A,B \}$, be the moment generating function of~$D$ conditional on the event~$E^{\alpha,\beta,C_2} \cap [S=C_1]$, and let~$\delta_{C_1,C_2}=1$ if $C_1=C_2$ and~$0$ otherwise. 
Then 
\begin{equation*}
M_{\pla}^{\alpha,\beta,C}(t) = 
\Big(\frac{(1-q)e^t}{1-qe^{2t}}\Big)^{\alpha+\beta}
\,
{\rm E}[e^{t(2R_{\pla}^{\alpha,\beta,C}-\delta_{B,C})}]
,\end{equation*}
for~$C\in\{ A,B \}$.
\end{theorem}
\vspace{1mm}

Quite remarkably, those moment generating functions (hence also all resulting moments) depend on~$(p_a,p_b)$ through~$q=(1-p_a)(1-p_b)$ only. Taking first and second derivatives with respect to $t$ in the above expressions and setting~$t=0$ then directly yields the following closed form expressions for the expected values~$e_{C_1}^{\alpha,\beta,C_2}$ from~(\ref{defexp}) and for the corresponding variances
$$
v_{C_1}^{\alpha,\beta,C_2}
{:=
\rm Var}[D\, |\, E^{\alpha,\beta,C_2}, S=C_1],
\qquad C_1,C_2\in\{ A,B \}.
$$

\begin{corollary} \label{expectscore} 
For $C\in\{A,B\}$, we have
(i)
$
\textstyle
e_{\pla}^{\alpha,\beta,C}
=
(\alpha+\beta)\,
\frac{1+q}{1-q}
-
\delta_{B,C}
+
2\,
{\rm E}[R_{\pla}^{\alpha,\beta,C}]
$
and
(ii) 
$
\textstyle
v_{\pla}^{\alpha,\beta,C}
=
4(\alpha+\beta)\,
\frac{q}{(1-q)^2}
+
4\,
{\rm Var}[R_{\pla}^{\alpha,\beta,C}]
$. 
Moreover, 
(iii) 
$e_{\pla}^{\alpha,\beta,C}$ is strictly monotone increasing in~$q$.
\end{corollary}
\vspace{2mm}

Clearly, Corollary~\ref{expectscore} confirms Simmons' result that the expected number of rallies  in an $\pla$-game won by $\pla$ on the score~$(n,k)$ is $e_{\pla}^{n,k,A}=
n\,
\frac{1+q}{1-q}
$ {for~$k=0$}. More interestingly, it also shows that the {exact value} for any~$k>0$ is given by
\begin{equation} \label{a}
\textstyle
e_{\pla}^{n,k,A}
=
(n+k)\,
\frac{1+q}{1-q}
+
2\sum_{r={1}}^{k}\, r  \,  W^{n,k,\pla}_\pla(q,r)  , 
\quad k=1,\ldots,n-1.
\end{equation}
Note that this is compatible with Simmons' result in~(\ref{simmons}) since the second term in the right-hand side of~(\ref{a}) is a weighted mean of $2r$, $r=1,\ldots, k$. Similarly, the expected number of rallies  in an $\pla$-game won by $\plb$ on the score~$(k,n)$, $k=0,1,\ldots,n-1$, is
$
\textstyle
e_{\pla}^{k,n,B}
=
(n+k)\,\frac{1+q}{1-q}-1
+2
\sum_{r={1}}^{k+1}\, r  \,  W^{k,n,\plb}_\pla(q,r)
$.

{The} expectation and variance of~$D$, in a $C$-game won by~$A$, are  {then} given by
\begin{equation} \label{eqE1}
\Bigg\{\
\begin{array}{l}
e_{C}^{A}
:=
{\rm E}[ D | E^{A}, S=C ]
=\frac{1}{p_C^A}\, \sum_{k=0}^{n-1} p_{C}^{n,k,\pla} e_C^{n,k,A}
\\[3mm]
v_{C}^{A}
:=
{\rm Var}[ D | E^{A}, S=C ]
=
\Big[
\frac{1}{p_C^A}\ \sum_{k=0}^{n-1} p_{C}^{n,k,\pla} (v_C^{n,k,A}+(e_C^{n,k,A})^2) 
\Big]
-  (e_{C}^{A})^2,
\end{array}
\end{equation}
while, in a $C$-game unconditional on the winner, they are given by
\begin{equation} \label{eqE2}
\Bigg\{\
\begin{array}{l}
e_C
:={\rm E}[D | S=C]
=
p_C^{\pla} e_{C}^{\pla} + p_C^{\plb} e_{C}^{B}
,
\\[2mm]
v_C:={\rm Var}[D | S=C]
=
  ( v_C^A + (e_C^A)^2) p_C^A+ 
  ( v_C^B + (e_C^B)^2)  p_C^B 
-  (e_{C})^2.
\end{array}
\end{equation}
Finally, unconditional on the server, this yields
\begin{equation}\label{eqE3}
\left\{\
\begin{array}{l}
e^A
:={\rm E}[D | E^A]
=
e_A^A s_a
+e_B^A s_b
,
\quad 
e
:={\rm E}[D]
= 
e_A s_a + 
e_B s_b,
\\[2mm]
v^A
:={\rm Var}[D | E^A]
=
(v_A^A+(e_A^A)^2) s_a
+(v_B^A+(e_B^A)^2) s_b
- (e^A)^2
,
\\[2mm]
v
:=
{\rm Var}[D]
=
(v_A+ e_A^2) s_a + 
(v_B+ e_B^2) s_b  - e^2
.
\end{array}
\right.
\end{equation}

Figures~\ref{fig1}(c)-(f) plot, for $n=15$, $e_{\pla}^{n,k,A}$, $e_{\pla}^{k,n,B}$,~$(v_{\pla}^{n,k,A})^{1/2}$, and $(v_{\pla}^{k,n,B})^{1/2}$ versus~$k$ for~$(p_a,p_b)=(.7,.5)$, $(.6,.5)$, $(.5,.5)$, and $(.4,.5)$, and report the corresponding numerical values of~$e_{\pla}^A$, $e_{\pla}^B$, $e_{\pla}$, $(v_{\pla}^A)^{1/2}$, $(v_{\pla}^B)^{1/2}$,  and $v_{\pla}^{1/2}$.
All expectation and standard deviation curves appear to be  strictly monotone increasing functions of the number $(n+k)$ of points scored, which was maybe expected. More surprising is the fact that---if one discards very  small values of~$k$---these curves are also roughly linear. Clearly, Simmons' lower and upper bounds~(\ref{simmons}), which are plotted versus $k$ in Figure~\ref{fig1}(c), only provide poor approximations of the exact expected values, particularly so for large~$k$. 

The dependence on~$(p_a,p_b)$ may be more interesting than that on~$k$. Note that, for each~$k$, $e_{\pla}^{n,k,A}$ and $e_{\pla}^{k,n,B}$ (hence also, $e_{\pla}^A$, $e_{\pla}^B$, and $e_{\pla}$) are decreasing functions of~$p_a$, which confirms Corollary~\ref{expectscore}(iii). Similarly,  
all quantities related to standard deviations also seem to be decreasing functions of~$p_a$. Now, it is seen  that, as a function of $p_a$, the expectation $e_\pla^\pla$ is more spread out than $e_\pla^\plb$. Indeed, the former ranges from $32.95$ ($p_a=.7$) to $56.30$ ($p_a=.4$), whereas the latter ranges from $41.95$ to $51.43$. On the contrary, the standard deviation of $D$ is more concentrated in an $\pla$-game won by $\pla$  (where it ranges from 8.34  ($p_a=.7$)  to 10.90  ($p_a=.4$))   than in an $\pla$-game won by $\plb$  (where it ranges from 7.36 to~11.44).  This phenomenon will appear even more clearly in Figure~\ref{fig3} below, where the same values of $(p_a, p_b)$ are considered. 
Note that the values of $e_{\pla}^A$, $e_{\pla}^B$, and $e_A$ are totally in line with  the score distribution and the expected values of~$D$ for each scores. For instance, the value~$e_{\pla}^B=41.95$ for~$p_a=.7$ translates the fact that when $B$ wins such an $A$-game, it is very likely (see Figure~\ref{fig1}(b)) that he/she will do so on a score that is quite tight, resulting on a large expected value for~$D$ (whereas, a priori, the values of $e_{\pla}^{k,n,B}$ range from 47.82 to 21.29 when $k$ goes from 14 to 0). The dependence of the expectation and standard deviation of~$D$ on rally-winning probabilities will further be investigated in Section~\ref{amerisec} for the no-server model when comparing the side-out scoring system with its rally-point counterpart.

Finally, in the case  $p_a=p_b=.5$, the fact  {that} $A$ is the first server in the game again brings some asymmetry  in the expected values and standard deviations of~$D$; in particular, this serve advantage alone is responsible for the fact that~$48.31=e_A^A<e_A^B=49.17$, and, maybe more mysteriously, that~$10.23=(v_{\pla}^A)^{1/2}>(v_{\pla}^B)^{1/2}=9.95$. 

\vspace{1mm}

\subsection{Distribution.}
\label{distriduration}
\vspace{2mm}

The  moment generating functions given in Theorem \ref{momentgenerating} allow, through a suitable change of variables,  for  obtaining the corresponding probability generating functions. These can in turn  be rewritten as power series whose coefficients yield the distribution of~$D$ conditional on the event~$E^{\alpha,\beta,C}\cap [S=A]$ (see the Appendix for the proof).

\begin{theorem} \label{probabilitygenerating}
Let $z\mapsto G_{C_1}^{\alpha,\beta,C_2}(z)={\rm E}[ z^D\, |\, E^{\alpha,\beta,C_2}, S=C_1 ]$, $C_1,C_2\in\{ A,B \}$, be the probability generating function of~$D$ conditional on the event~$E^{\alpha,\beta,C_2}\cap [S=C_1]$. Then, for~$C \in\{ A,B \}$, 
$$
\textstyle
G_\pla^{\alpha, \beta, C}(z) 
= 
\frac{p_a^{\alpha} p_b^{\beta} q_a^{\delta_{B,C}}}{p_{\pla}^{\alpha,\beta,C}}\,
\sum_{j=0}^\infty \
 q^j H_{\pla}^{\alpha,\beta,C}(j)
\,
z^{\alpha+\beta+2j+\delta_{B,C}},
$$ 
 where, writing $m^+:=\max(m,0)$, we let
$$
\textstyle
H_{\pla}^{\alpha,\beta,\pla}(j) :=\textstyle \sum_{l=(j-\gamma_1)^+}^{j} {{\alpha+\beta+l-1}\choose{l}} 
\binom{\alpha}{j-l} \binom{\beta-1}{j-l-1}
$$
and
$$
\textstyle
H_{\pla}^{\alpha,\beta,\plb}(j) := \sum_{l=(j-\gamma_2 )^+}^{j} {{\alpha+\beta+l-1}\choose{l}} 
\binom{\alpha}{j-l} \binom{\beta-1}{j-l}.$$
\end{theorem}
\vspace{2mm}

This result gives the probability distribution of $D$, conditional on~$E^{\alpha,\beta,C}\cap\, [S=A]$, for~$C\in\{A,B\}$. Note that, as expected, we have ${\rm P}[D=d\, | \,E^{\alpha,\beta,A}, S=A]=0={\rm P}[D=d+1 \,|\, E^{\alpha,\beta,B}, S=A]$ for all $d < \alpha+\beta$.  Moreover,  for all nonnegative integer~$j$,  ${\rm P}[D=\alpha+\beta+2j+1 \; | \; E^{\alpha,\beta, \pla}, S=\pla]=0 = {\rm P}[D=\alpha+\beta+2j \; | \; E^{\alpha,\beta, \plb}, S=\pla]$. In the sequel, we refer to this  as the \emph{server-effect}.

Theorem \ref{probabilitygenerating} of course allows for investigating the shape of the distribution of~$D$ above all scores, and not only, as in Figures~\ref{fig1}(c)-(f), its expectation and standard deviation. This is what is done in  Figure~\ref {fig2}, which plots, as a function of the score, quantiles of order~$\alpha=.01$, .05, .25, .5, .75, .95, and~.99 for~$(p_a,p_b)=(.6,.5)$. For each~$\alpha$, two types of quantiles are reported, namely (i) the standard quantile~$q_\alpha:=\inf \{ d : {\rm P}[D\leq d\,|\, E^{\alpha,\beta,C}, S=A] \geq \alpha \}$ and (ii) an interpolated quantile, for which the interpolation is conducted linearly over the set~$(d,d+2)$ containing the expected quantile (here, we avoid interpolating over ($d,d+1$) because of the above server-effect, which implies that either $d$ or $d+1$ does not bear any probability mass). One of the most prominent features of Figure~\ref{fig2} is the wiggliness of the standard quantile curves, which is directly associated with the server-effect. It should be noted that the expectation curves (which are the same as in Figures~\ref{fig1}(c)-(d)) stand slightly above the median curves, which possibly indicates that, above each score, the conditional distribution of~$D$ is somewhat asymmetric to the right. This  {(light) asymmetry} is confirmed by the other quantiles curves.

Now, the probability distribution of~$D$ in an $\pla$-game,  unconditional on the score, is of course derived trivially from its conditional version obtained above and the score distribution of Section~\ref{probagame}. The general form of this distribution is somewhat obscure (and will not be explicitly given here), but it yields easily interpretable expressions for small values of $d$. For instance, one obtains
\begin{gather*}
{\rm P}[D=n | S=A] =p_a^n, 
\\[1mm] 
{\rm P}[D=n+1 | S=A] =  q_a p_b^n, 
\\[1mm] 
{\rm P}[D=n+2 | S=A] = nqp_a^n+p_a q_ap_b^n,\ldots
\end{gather*}
Finally,  the unconditional distribution of~$D$ is simply obtained through
$
{\rm P}[D=k]
=
{\rm P}[D=k | S=A]
s_a
+
{\rm P}[D=k | S=B]
s_b$,
$k\geq0$, where one computes the distribution for a $\plb$-game by inverting $p_a$ and $p_b$ in the distribution for an $\pla$-game.  
 
Figure~\ref{fig3} shows that  {there are a number} of remarkable aspects to these distributions.  First note the influence of the above mentioned server-effect, which causes the wiggliness visible in most curves there. Also note that the distributions  in Figure~\ref{fig3}(c) are much less wiggling than the corresponding curves in Figures~\ref{fig3}(a)-(b). As it turns out, this wiggliness is present, albeit more or less markedly, at all stages (that is, not only to the right of the mode)   for every choice of $(p_a, p_b)$. 
Most importantly, despite their irregular aspect, all curves are essentially unimodal, as expected.
Now, consider  the dependence  on~$p_a$ of the position and spread of these curves. One sees that  while their spread clearly increases much more rapidly with~$p_a$ in Figure~\ref{fig3}(b) than  in Figure~\ref{fig3}(a),  the opposite can be said for their mode. This is easily understood in view of the corresponding means and variances, which are recalled in the legend boxes (and coincide with those from  Figure~\ref{fig1}). As for the curves in Figure~\ref{fig3}(c), they are obtained by averaging the corresponding  curves in Figure~\ref{fig3}(a) and Figure~\ref{fig3}(b) with weights 
$p_A^A$
and 
$p_A^B=1-p_A^A$,
respectively. Taking into account the values of these probabilities  explains why  the curves with $p_a=.7$ and $p_a=.6$     are essentially the corresponding curves  in Figure~\ref{fig3}(a), whereas that with $p_a=.4$ is closer to the corresponding curve in Figure~\ref{fig3}(b). 

\section{Comparison with the rally-point scoring system.} \label{amerisec}
\vspace{2mm}

One of the main motivations for this work was to compare more deeply the side-out scoring system considered in Sections~\ref{probagame} and \ref{duration} with the \emph{rally-point scoring system}. 
As mentioned in the Introduction, many sports recently switched (e.g., badminton, volley-ball)---or are in the process of switching (e.g., squash)---from the side-out scoring system to its rally-point counterpart{, whereas others (e.g., racquetball) so far are sticking} to the side-out scoring system. It is therefore natural to investigate the implications of the transition to the rally-point system. 

The literature, however, has focused on the impact of the scoring system on the \emph{outcome} of the game---studied by comparing the game-winning probabilities under both scoring systems; see, e.g., Simmons (1989). This is all the more surprising since there have been, in the sport community, much debate and questions about how much the duration of the game is affected by the scoring system. Moreover, it is usually reported that the main motivation for turning to the rally-point system is  to regulate the playing time (that is,  to make the length of the match more predictable), which is of primary importance for television, for instance. Whether the transition to the rally-point system has indeed served that goal, and, if it has, to what extent, are questions that have not been considered in the literature, and were at best addressed on empirical grounds only {(by international sport federations).} 

In this section, we will provide an in-depth comparison of the two scoring systems, both in terms of game-winning probabilities and in terms of durations, which will provide theoretical answers to the questions above. Again, this is made possible by our rally-level analysis of the game and the results of the previous sections on the distribution of the number of rallies  under the side-out scoring system. As we will discuss in Section~\ref{finalcom}, our results are potentially of high interest both for international federations and for local tournament organizers.

\subsection{Game-winning probabilities.}
\label{probacompare}
\vspace{2mm}

Although the game-winning probabilities for an $A$-game played under the rally-point system have already been obtained in the literature (see, e.g., Simmons~1989), we start by deriving them quickly, mainly for the sake of completeness, but also because they easily follow from the combinatorial methods used in the previous sections. First note that there cannot be \emph{exchanges} (in the sense of Definition~\ref{defone}) in the rally-point scoring system. We then denote by 
$
\bar{E}_{\pla}^{\alpha,\beta,C}(r)
$
($C\in \{A,B\}$) 
the event associated with a sequence of rallies  that, in the rally-point system, (i) gives raise to $\alpha$ points  scored by Player~A and  $\beta$ points scored by  Player~B, (ii) involves exactly~$r$ $\pla$-interruptions, and (iii) is such that Player~$C$ scores a point in the last rally; again, it is  assumed here that $\alpha>0$ (resp., $\beta>0$) if $C=A$ (resp., if~$C=B$). We write 
$$
\bar{p}_{C_1}^{\alpha,\beta,C_2}(r) := {\rm P}[ \bar{E}^{\alpha,\beta,C_2}(r) | S=C_1 ],
\qquad
C_1,C_2\in \{A,B\}.
$$
The following result then follows along the same lines as for Lemma~\ref{probascoreinterexch} and Theorem~\ref{probascore}.

\begin{theorem} \label{probascoreamer}
(i) With the notation above, $
\textstyle
\bar{p}_{\pla}^{\alpha,\beta,\pla}(r)
=
{\alpha \choose r} {\beta-1 \choose r-1} p_a^{\alpha-r} p_b^{\beta-r} (q_a q_b)^{r}
$,
 $r\in\{\gamma_0,\ldots, \gamma_1\}$,
and~$
\textstyle
\bar{p}_{\pla}^{\alpha,\beta,\plb}(r)=
{\alpha \choose r-1} {\beta-1 \choose r-1} p_a^{\alpha-r+1} p_b^{\beta-r} q_a (q_a q_b)^{r-1},
$
$r\in\{ 1,\ldots, \gamma_2+1 \}$.
(ii) Writing $\bar{p}_{\pla}^{\alpha,\beta,C}$ for the probability of the event~$\bar{E}_{\pla}^{\alpha,\beta,C}:=\cup_{r} \, \bar{E}_{\pla}^{\alpha,\beta,C}(r)$, we have~$
\bar{p}_{\pla}^{\alpha,\beta,\pla}
= 
p_a^{\alpha} p_b^{\beta} \sum_{r=\gamma_0}^{\gamma_1} \,
{\alpha \choose r} {\beta-1 \choose r-1}  (t_a t_b)^{r}
$
and
$
\bar{p}_{\pla}^{\alpha,\beta,\plb}
= 
p_a^{\alpha} p_b^{\beta-1}q_a \sum_{r=1}^{\gamma_2+1} 
{\alpha \choose r-1} {\beta-1 \choose r-1}  
\linebreak
(t_a t_b)^{r-1}
$,
where we let $t_a=q_a/p_a$ and $t_b=q_b/p_b$.
\end{theorem}
\vspace{0mm}

\begin{remark} \label{amernoserve} These  {expressions} further simplify in the no-server model $(p:=)p_a=1-p_b$. There we indeed have $t_b= t_a^{-1}$, so that  the above formulas yield  $\bar{p}_{\pla}^{\alpha,\beta,\pla}= {\alpha+\beta-1 \choose \beta} p^{\alpha} (1-p)^{\beta} $  and   $\bar{p}_{\pla}^{\alpha,\beta,\plb} = {\alpha+\beta-1 \choose \alpha}p^{\alpha} (1-p)^{\beta}.$   
\end{remark}

Of course, the resulting score distribution and game-winning probabilities for an $A$-game directly follow from Theorem~\ref{probascoreamer}. In accordance with the notation adopted for the side-out scoring system, we will write
\vspace{-2mm}
$$
\bar{p}_C^{\pla}
:={\rm P}[\bar{E}^A | S=C]
:={\rm P}[\cup_{k=0}^{n-1} E^{n,k,A} | S=C]
:=\sum_{k=0}^{n-1} \bar{p}_C^{n,k,\pla},
\qquad
\bar{p}_C^{\plb}:=1-\bar{p}_C^{\pla},
$$
$$
\bar{p}^{n,k,A}:={\rm P}[\bar{E}^{n,k,A}]=\bar{p}_A^{n,k,A} s_a+ \bar{p}_B^{n,k,A} s_b,
\bar{p}^{k,n,B}:={\rm P}[\bar{E}^{k,n,B}]=\bar{p}_A^{k,n,B} s_a+ \bar{p}_B^{k,n,B} s_b,
$$
and
$$
\bar{p}^{C}:={\rm P}[\bar{E}^C]=\bar{p}_A^C s_a+\bar{p}_B^Cs_b.
$$

Figures~\ref{fig4}(a)-(b) plot the same score distribution curves as in Figures~\ref{fig1}(a)-(b), respectively, but in the case of an $A$-game played under the  rally-point scoring system with $n=21$. Both pairs of plots look roughly similar, although extreme scores seem to be less likely in the rally-point scoring; this confirms the findings from Simmons (1989) according to which shutouts are less frequent under the rally-point scoring system. Note also that, unlike for the side-out scoring, the $(p_a,p_b)=(.5,.5)$ curve in Figure~\ref{fig4}(a) is the exact reverse image of the corresponding one in Figure~\ref{fig4}(b): for the rally-point scoring, Player~A does not get any advantage from serving first if~$(p_a,p_b)=(.5,.5)$, which is confirmed by the game-winning probabilities~$\bar{p}_A^A=\bar{p}_A^B=.5$. 

Again, the dependence of the game-winning probabilities   {on}~$(p_a,p_b)$ is of primary importance. We will investigate this dependence visually and compare it with the corresponding dependence for the side-out scoring system. To do so, we focus on the no-server version ($p=p_a=1-p_b$) of Badminton, where, as already mentioned, the side-out scoring system with~$n=15$ (men's singles) was recently replaced with the rally-point one characterized by~$n=21$. The results are reported in Figures~\ref{fig5}(a)-(b). Figure~\ref{fig5}(a) supports the claim---reported, e.g., in Simmons (1989) or Percy~(2009)---stating that, for any fixed~$p$, the scoring barely influences game-winning probabilities.  Now, while Figure~\ref{fig5}(b) shows that the probability that Player~$A$ wins an $A$-game is essentially the same for both scoring systems if he/she is the best player ($\bar{p}_A^A/p_A^A\in(.926,1)$ for $p\geq.5$, and $\bar{p}_A^A/p_A^A\in(.997,1)$ for~$p>.7$), it tells another story for $p<.5$: there,   the probability that~$A$ wins an $A$-game played under the rally-point system (i)  becomes relatively negligible for very small values of~$p$ (in the sense that $\bar{p}_A^A/p_A^A\to 0$ as~$p\to 0$) and (ii) can be up to~28 times larger than under the  side-out system (for values of~$p$ close to $.1$). Of course, one can say that (i) is irrelevant since it is associated with an event (namely, a victory of $A$) occurring with very small probability; (ii), however, constitutes an important difference between both scoring systems for values of~$p$ that are not so extreme. 

\subsection{Durations.}
\label{durationcompare}
\vspace{2mm}

In the rally-point system, the number of rallies  needed to achieve the event~$\bar{E}^{\alpha,\beta,C_2}\cap[S=C_1]$ is not random: with obvious notation, it is almost surely equal to $\bar{e}_{C_1}^{\alpha,\beta,C_2}=\alpha+\beta$, which explains why  {Figure~\ref{fig4}} does not contain the rally-point counterparts of Figure~\ref{fig1}(c)-(f). 
The various conditional and unconditional means and variances of the number of rallies  in the rally-point system (that is, the quantities $\bar{e}_{C}^{A}$, $\bar{v}_{C}^{A}$, $\bar{e}_C$, $\bar{v}_C$, $\bar{e}^A$, $\bar{v}^A$, $\bar{e}$, $\bar{v}$) can then be readily computed from the game-winning probabilities given in Theorem~\ref{probascoreamer}, in the exact same way as in~(\ref{eqE1})-(\ref{eqE3}) for the side-out scoring system. More generally, the corresponding distribution of the number of rallies  in a game trivially follows from the same game-winning probabilities.
 

Figures~\ref{fig5}(c)-(h) plot, as functions of $p=p_a=1-p_b$ (hence, in the no-server model), expected values and standard deviations of the numbers of rallies  needed to complete (i) $A$-games played under the side-out system with $n=15$ and (ii)  $A$-games played under the rally-point system with $n=21$. Clearly, those plots allow for an in-depth original comparison of both scoring systems. Let us first focus on durations unconditional on the winner of the game. Figure~\ref{fig5}(c) shows that (i)  games played under the side-out system will last longer than those played under the rally-point one for players of roughly the same level (which was expected since the side-out system will then lead to many exchanges), whereas (ii) the opposite is true when one player is much stronger 
  (which is explained by the fact that shutouts require more rallies  in the rally-point scoring considered than in the side-out one). Maybe less expected is the fact (Figure~\ref{fig5}(f)) that the standard deviation of~$D$ is, uniformly in~$p\in (0,1)$, smaller for the rally-point scoring system than for the side-out system, which shows that the transition to the rally-point system indeed makes the length of the match more predictable.  The twin-peak shape of both standard deviation curves is even more surprising. Finally, note that, while the rally-point curves  in Figures~\ref{fig5}(c) and (f) are symmetric about~$p=.5$, the side-out curves are not, which is due to the server-effect. This materializes into the limits of $e_A$ given by $16$ and $15$ as $p\to 0$ and $p\to 1$, respectively (which was expected: if Player~$B$ wins each rally with probability one, he/she will indeed need 16 rallies  to win an $A$-game, since he/she has to regain the right to serve before scoring his/her first point),  but also translates into (i) the fact that the mode of the side-out curve in Figure~\ref{fig5}(c) is not exactly located in $p=.5$ and (ii) the slightly different heights of the two local (side-out) maxima in Figure~\ref{fig5}(f). 

We then turn to durations conditional on the winner of the game, whose expected values and standard deviations are reported in Figures~\ref{fig5}(d), (e), (g), and (h). These figures look most interesting and reveal important differences between both scoring systems. Even the general shape of the curves there are of a different nature for both scorings; for instance, the rally-point curves in Figures~\ref{fig5}(d)-(e) are monotonic, while the side-out ones are unimodal. Similarly, in Figure~\ref{fig5}(g), the rally-point curve is unimodal, whereas the side-out curve exhibits a most unexpected bimodal shape. It is also interesting to look at limits as $p\to 0$ or $p\to 1$ in those four subfigures; these limits, which are derived in Appendix~\ref{limitcalc}, are plotted as short horizontal lines.
Consider first limits above events occurring with probability one, that is, limits as $p\to 1$ in Figures~\ref{fig5}(d), (g) and limits $p\to 0$ in Figures~\ref{fig5}(e), (h). The resulting limits are totally in line with the intuition: the four conditional standard deviations go to zero, which implies that the limiting conditional distribution of~$D$ simply is almost surely equal to the corresponding limiting (conditional) expectations. The latter themselves assume very natural values: for instance, for the same reason as above, $e_A^B$ converges to~$16$, which is therefore the limit of~$D$ in probability.

Much more surprising is what happens for limits above events occurring with probability zero, that is, limits as $p\to 0$ in Figures~\ref{fig5}(d), (g) and limits as $p\to 1$ in Figures~\ref{fig5}(e), (h). Focussing first on the side-out scoring system, it is seen that a (miraculous) victory of $A$ will require, in the limit, almost surely~$D=15$ points, while the limiting conditional distribution of~$D$ for victories of~$B$ is non-degenerate. The latter distribution is shown (see Appendix~\ref{limitcalc}) to be \emph{uniform over~$\{n+1,n+2,\ldots,2n\}$} (hence is stochastically bounded!), which is compatible with the values~$n+1+(n-1)/2(\approx 3n/2)$ and~$(n-1)^2/12$ for the limiting expectation and variance, respectively. It should be noted here that this huge difference between those two limiting conditional distributions of~$D$ is entirely due to the server-effect. In the absence of the server-effect, the subfigures (e) and (h) should indeed be the exact reverse image of the subfigures (d) and (g), respectively.  Similarly, the bimodality of the side-out curve in Figure~\ref{fig5}(g) is also due to the server-effect. We then consider the rally-point scoring, which is  not affected by the server-effect, so that it is sufficient to consider at the limits as~$p\to 0$ in Figures~\ref{fig5}(d), (g). There, one also gets a non-degenerate limiting conditional distribution for~$D$, with expectation~$2n^2/(n+1)(\approx 2n)$ and variance~$2n^2(n-1)/[(n+1)^2(n+2)](\approx 2)$.  

\section{Simulations.}\label{simu}

We performed several Monte Carlo simulations, one for each figure considered so far (except Figure~\ref{fig2}, as it already contains many theoretical curves). 
To describe the general procedure, we focus on the Monte Carlo experiment associated with the side-out scoring system in Figure~\ref{fig5} (results for the rally-point scoring system there or for the other figures are obtained similarly). For each of the $1,999$ values of~$p$ considered in Figure~\ref{fig5}, the corresponding values of~$p_A^A(p)$, $e_A(p)$, $v_A(p)$, $e_A^C(p)$,  $v_A^C(p)$, $C\in\{A,B\}$,  were estimated on the basis of $J=200$ independent replications of an $A$-game played under the side-out scoring system with~$p_a=1-p_b=p$. Of course, for each fixed~$p$, the game-winning probability~$p_A^C(p)$ is  simply estimated by the proportion of games won by~$C$ in the  $J$ corresponding  $A$-games:
$$
\hat{p}_A^C(p)
:=\frac{J^C}{J}
:=\frac{1}{J} \sum_{j=1}^J I_{j}^C,
$$ 
where $I_j^C$, $j=1,\ldots,J$, is equal to one (resp., zero) if Player~$C$ won (resp., lost) the $j$th game. The corresponding estimates  for~$e_A(p)$, $v_A(p)$, $e_A^C(p)$, and $v_A^C(p)$ then are  given by
$$
\hat{e}_A(p)
:=
\frac{1}{J} \sum_{j=1}^J  d_j  
\quad
\quad
\hat{v}_A(p)
:=
\frac{1}{J} \sum_{j=1}^J  \big( d_j - \hat{e}_A(p)\big)^2,
$$
\begin{equation} \label{condexv}
\hat{e}_A^C(p)
:=
\frac{1}{J^C} \sum_{j=1}^J  d_j I_{j}^C ,
\quad
\textrm{and}
\quad
\hat{v}_A^C(p)
:=
\frac{1}{J^C} \sum_{j=1}^J  \big( d_j - \hat{e}_A^C(p)\big)^2 I_{j}^C,
\end{equation}
where~$d_j$, $j=1,\ldots,J$, is the total number of rallies  in the $j$th game. These estimates are plotted in thin blue lines in Figure~\ref{fig5}. Clearly, these simulations validate our theoretical results in Figures~\ref{fig5}(a), (c), and~(f). To describe what happens in the other plots, consider, e.g., Figure~\ref{fig5}(g). There, it appears that the theoretical results are confirmed for large values of~$p$ only. However, this is simply explained by the fact that for small values of~$p$, the denominator of~$\hat{v}_A^A(p)$ (see~(\ref{condexv})) is very small. Actually, among the $542\times 200$ $A$-games associated with the 542 values of $p\leq .2710$, not a single one here led to a victory of~$A$, so that the corresponding  estimates~$\hat{v}_A^A(p)$ are not even defined. Of course, values of~$p$ slightly larger than $.2710$  still give raise to a small number of victories of~$A$, so that the corresponding estimates~$\hat{v}_A^A(p)$ are highly unreliable. The situation of course improves substantially as~$p$ increases, as it can be seen in Figure~\ref{fig5}(g). Figures~\ref{fig5}(b), (d), (e), and~(h) can be interpreted exactly in the same way.

This underlines the fact that expectations and variances conditional on events with small probabilities are of course extremely difficult---if not impossible---to estimate. To quantify this, let us focus again on Figure~\ref{fig5}(g), and consider the local maximum on the left of the plot, which is (on the grid of values of~$p$ at hand) located in~$p_0:=.0085$. The probability~$p_A^A(p_0)$ of a victory of~$A$ in an $A$-game played under the side-out scoring system with~$p=p_0$ is about $3.5\times10^{-31}$. Estimating $v_A^A(p_0)$ with the same accuracy as that  achieved for, e.g., $v_A^A(.5)$ in Figure~\ref{fig5}(g) would then require a number of replications of (fixed $p_0$) $A$-games that is about~$200\times p_A^A(.5)/p_A^A(p_0)\approx 3\times10^{32}$. Assuming that $10^6$ replications can be performed in a second by a super computer (which is overly optimistic), this estimation of $v_A^A(p_0)$ would still require not less than $9.5\times 10^{18}$ years! This means that it is indeed impossible to estimate in a reliable way the conditional variance curve for~$p$ close to~$p_0$, hence that there is no way to empirically capture the convergence of $v_A^A(p)$ to~$0$ as $p\to 0$. Without our theoretical analysis, there is therefore no hope to learn about the degeneracy (resp., non-degeneracy) of the limiting distribution of~$D$ conditional on a victory of~$A$ as~$p\to 0$ (resp.,  conditional on a victory of~$B$ as~$p\to 1$).

We will not comment in detail the Monte Carlo results associated with the other figures. We just report that they again confirm our theoretical findings, whenever possible, that is, whenever they are not associated with conditional results above events with small probabilities. 
\vspace{3mm}

\section{Conclusion and final comments.}
\label{finalcom}

This paper provides a \emph{complete} rally-level probabilistic description for games played under the side-out scoring system. It complements the previous main contributions by Phillips~(1978), Strauss and Arnold~(1987), and Simmons (1989) by adding to the well-known game-winning probabilities an exhaustive knowledge of the random duration of the game. 
This brings  a much better understanding of the underlying process as a whole, as is demonstrated in Sections~\ref{probagame} to \ref{amerisec}.

In this final section, we will mainly focus on the practical implications of our findings. For this, we may restrict to~$(p_a,p_b)\in[.4,.6]\times[.4,.6]$, say, since players tend to be grouped according to strength. For such values of the rally-winning probabilities, our results
show that the recent transition---in mens' singles' Badminton---from the $n=15$ side-out scoring system to the $n=21$ rally-point one strongly affected the properties of the game. They indeed 
indicate that (i) games played under the rally-point scoring system are much shorter than those played according to the side-out one, and that (ii) the uncertainty in the duration of the match is significantly reduced. 
Our results allow to quantify both effects. On the other hand, they show that game-winning probabilities are essentially the same for both scoring systems. It is then tempting to conclude (as in Simmons~(1989) and Percy~(2009)) that the outcomes of the matches are  barely influenced by the scoring system adopted.  {While this is strictly valid in the model, it is highly disputable under possible violations of the model (stating, e.g., that players might get tired at different speeds) which, given the reduced duration of the game emphasized by our analysis, may appear quite relevant.}


In practice, the results of this paper can be useful to many actors of the sport community. For the international sport federations playing with the idea of replacing the side-out scoring system with the rally-point one, our results could be used to tune~$n$ (i.e., the number of points to be scored to win a rally-point game) according to their wishes. For the sake of illustration, consider again the transition performed by the International Badminton Federation (IBF). Presumably, their objective was (i) to make the duration of the game more predictable and (ii) to ensure that  the outcome of the matches would change as little as possible. If this was indeed their objective, then our results show that it has been partially achieved.  {However, it is now easy to see that other choices of~$n$ would have been even better in that respect, the choice of $n=27$ (see Figure~\ref{fig6}(d) and (b)), being optimal. Moreover,    this last choice would have affected  the duration of the game much less than $n=21$ (see Figure~\ref{fig6}(c)), and thus would have made the outcome of the matches more robust to possible violations of the model.}


For organizers of local tournaments played under the side-out scoring system, our results can be used to control, for any fixed number of planned matches, the time required to complete their events. Such a control over this random time, at any fixed tolerance level, can indeed be achieved in a quite direct way from our results on the duration of a game played under the side-out scoring system. Organizers can then deduce, at the corresponding tolerance level, the number of matches---hence the number of players---their events will be able to host. This of course concerns the sports that are still using this scoring system, such as racquetball and squash (for the latter, only in countries currently using the so-called English scoring system).

Finally, our results also open the way to more efficient estimation of the rally-winning probabilities~$(p_a,p_b)$ in the side-out scoring system. Consequently, they potentially lead to more accurate \emph{ranking} procedures (based  on Bradley-Terry paired comparison methods), which is of course of high interest to national and international federations still supporting that scoring system. A full discussion of this is beyond the scope of this paper (and is actually the topic of Paindaveine and Swan~2009), and we only briefly describe the main idea here. 
Essentially, the results of the present work, in a point estimation framework, enable us to perform maximum likelihood estimation of ($p_a,p_b$) based on \emph{game scores and durations}. It is natural to wonder how much the resulting estimators would improve on the purely score-based maximum likelihood estimators proposed by  Strauss and Arnold~(1987). It turns out that the improvement is very important. First of all, unlike the purely score-based estimators, which require numerical optimization techniques, the score-and-duration-based ones happen to allow for elegant closed form expressions. Second, they can be shown to enjoy strong finite-sample optimality properties. Last but not least, they are much more accurate than their Strauss and Arnold (1987) competitors, especially so for small numbers of observations (i.e., games). This is illustrated in Figure~\ref{fig7}, where it can be seen that even for as little as $m=2$ games, the score-and-duration-based maximum likelihood estimators outperform their competitors both in terms of bias and variability.  Clearly, this will have priceless practical  consequences for the resulting ranking procedures, since any fixed pair of players typically do not meet more than once or twice in the period (usually one year) on which this ranking is to be computed.

\appendix
\section{Appendix.}

\subsection{Proof of Lemma~\ref{probascoreinterexch} and Theorem~\ref{probascore}.}
\vspace{1mm}

In the Appendix, we simply write \emph{interruptions} for \emph{$\pla$-interruptions}.
\vspace{2mm}

\noindent \textbf{Proof of Lemma~\ref{probascoreinterexch}.}
Clearly, $p_{\pla}^{\alpha,\beta,\pla}(r,j)=K_{r,j}\,p_a^{\alpha} p_b^{\beta} (q_a q_b)^{r+j}$, where $K_{r,j}$ is the number of ways of setting $r$ interruptions and $j$ exchanges in the sequence of rallies  achieving  the event under consideration. Regarding interruptions, we argue as in~Hsi and Burich~(1971), and say those $r$ interruptions should be put into the $\alpha$ possible spots (remember the last point should be won by $\pla$), while the $\beta$ points scored by~$\plb$ should be distributed among those $r$ interruptions---with at least one point scored by~$\plb$ in each interruption (so that there may be at most $r=\min(\alpha,\beta)$ interruptions). There are exactly ${\alpha \choose r} {\beta-1 \choose r-1}$ ways to achieve this.
 As for the $j$ exchanges, they may occur at any time and thus  there are as many ways of placing $j$ interruptions as there are distributions of $j$ indistinguishable balls into $\alpha+\beta$ urns, i.e. $\alpha+\beta-1\choose j$. 
Summing up, we have proved that 
$$
p_{\pla}^{\alpha,\beta,\pla}(r,j)
= {{\alpha+\beta-1}\choose j }{\alpha \choose r} {\beta-1 \choose r-1} p_a^{\alpha} p_b^{\beta} (q_a q_b)^{r+j},
$$
with $r=\mino{\beta}{1},\ldots, \mino{\alpha}{\beta}$, $j\in\Nn$.

As for $p_{\pla}^{\alpha,\beta,\plb}(r,j)$, this probability is clearly of the form $L_{r,j}\,p_a^{\alpha} p_b^{\beta} q_a(q_a q_b)^{r+j-1}$. In this case, there are $\alpha+1$ possible spots for the $r$ interruptions. But since $\plb$ scores the last point, the sequence of rallies  should end with an interruption. There are therefore ${\alpha \choose r-1}$ ways to insert the interruptions. Each interruption contains at least one point for $\plb$, so that $r\leq\min(\alpha+1,\beta)$. The result follows by noting that there are \vspace{-1mm} ${\beta-1 \choose r-1}$ ways of distributing the $\beta$ points scored by $\plb$ into those $r$ interruptions, and by dealing with exchanges as for $p_{\pla}^{\alpha,\beta,\pla}(r,j)$.
\cqfd
\vspace{3mm}

\noindent \textbf{Proof of Theorem~\ref{probascore}.}
The result directly follows from Lemma~\ref{probascoreinterexch} by writing $p_{\pla}^{\alpha,\beta,\pla}=\sum_{r,j} p_{\pla}^{\alpha,\beta,\pla}(r,j)$ and $p_{\pla}^{\alpha,\beta,\plb}=\sum_{r,j} p_{\pla}^{\alpha,\beta,\plb}(r,j)$ (where the sums are over all possible values of $r$ and~$j$ in each case), and by using the equality
$\sum_{j=0}^{\infty}  {{m+j-1}\choose j } z^{j}=(1-z)^{-m}$ for any $z\in[0,1)$.
\cqfd
\vspace{3mm}

\subsection{Proof of Theorems~\ref{momentgenerating} and~\ref{probabilitygenerating} and of~Corollary~\ref{expectscore}.}
\vspace{2mm}

\noindent \textbf{Proof of Theorem~\ref{momentgenerating}.}
First note that if $\pla$ scores the last point in an $\pla$-game in which the score is~$\alpha$ to $\beta$ after $j$ exchanges ($j \in \{0, 1, \ldots \}$)  and $r$ interruptions  ($r \in \{\gamma_0, \ldots, \gamma_1\}$), then there have been~$ \alpha + \beta +2(r+j)$ rallies . Conditioning on the number of interruptions  and exchanges therefore yields 
\[M_{\pla}^{\alpha,\beta,\pla}(t)=(p_\pla^{\alpha, \beta, \pla})^{-1}{\sum_j\sum_r}e^{t( \alpha + \beta +2(r+j))}p_\pla^{\alpha, \beta, \pla}(r,j)\] (where the sums are over all possible values of $r$ and~$j$ in each case) and thus, from Lemma~\ref{probascoreinterexch} and Theorem~\ref{probascore}, 
$$\begin{array}{rcl} 
M_{\pla}^{\alpha,\beta,\pla}(t) & = & \dfrac{ e^{t(\alpha+\beta)} 
\sum_j\left(e^{2t}q\right)^j\binom{\alpha+\beta+j-1}{j}\sum_re^{2tr}\binom{\alpha}{r}\binom{\beta-1}{r-1}q^r}{(1-q)^{-(\alpha+\beta)}\sum_r\binom{\alpha}{r}\binom{\beta-1}{r-1}q^r}\\
& & \\
& = & \left((1-q)e^t\right)^{\alpha+\beta}\left(\sum_j\left(e^{2t}q\right)^j\binom{\alpha+\beta+j-1}{j}\right)\left(\sum_re^{2tr}W_\pla^{\alpha, \beta, \pla}(q, r)\right). \end{array}$$
The first  claim of Theorem~\ref{momentgenerating} follows.  

For the second claim, it suffices to note that  if $\plb$ scores the last point in an $\pla$-game in which the score is of $\alpha$ to $\beta$ after $j$ exchanges ($j \in \{0, 1, \ldots \}$)  and $r$ interruptions  ($r \in \{1, \ldots,\gamma_2+1\}$),  then  the number of rallies  equals $\alpha+\beta+2(r-1+j)+1$; the computations above  then hold with only minor changes. 
\cqfd
\vspace{3mm}

\noindent \textbf{Proof of Corollary~\ref{expectscore}.} 
Taking first and second derivatives of the moment generating functions yields the expectations and variances given in Corollary~\ref{expectscore}. Moreover it   can easily be seen that derivatives of the expected values with respect to $q$ are positive by using the Cauchy-Schwarz inequality, and thus the latter are strictly monotone increasing in $q$. 
 \cqfd
\vspace{3mm}

\noindent \textbf{Proof of Theorem~\ref{probabilitygenerating}.} The change of variables $z = e^t$ in the moment generating functions given in Theorem~\ref{momentgenerating} immediately yields the probability generating functions. If $\beta=0$, the latter is already in the form of an infinite series
$G_\pla^{\alpha, 0, \pla}(z) = \sum_{j=0}^\infty(1-q)^\alpha q^j\binom{\alpha+j-1}{j}z^{\alpha+2j}.$
If $\beta>0$,  we have
$$G_\pla^{\alpha, \beta, \pla}(z) = (1-q)^{\alpha+\beta}z^{\alpha+\beta}\sum_{j=0}^\infty K_j z^{2j}\sum_{r=1}^{\gamma_1} W_r z^{2r},$$
where  $K_j = q^j\binom{\alpha+\beta+j-1}{j}$ and $W_r = W^{\alpha,\beta,\pla}_\pla(q,r)$.
This double sum satisfies
$$\sum_{j=0}^\infty K_j \sum_{r=1}^{\gamma_1} W_r z^{2(j+r)} = \sum_{j=1}^{\gamma_1} z^{2j}\left( \sum_{l=0}^{j-1}K_lW_{j-l}\right)+\sum_{j=\gamma_1+1}^\infty z^{2j}\left(\sum_{l=j-\gamma_1}^{j-1}K_lW_{j-l}\right).$$
The same arguments are readily adapted to $G_\pla^{\alpha, \beta, \plb}(z)$, and  Theorem~\ref{probabilitygenerating} follows. \cqfd

\subsection{The distribution of the number of rallies , in the no-server model, for extreme values of the rally-winning probabilities.}
\label{limitcalc}
\vspace{2mm}

As announced in Section~\ref{durationcompare}, we determine here the limiting behavior of the number of rallies~$D$, in the no-server model, for~$p\to0$ and~$p\to1$, conditional on the winner of the $A$-game considered. We start with the limit under almost sure events, that is, limits as~$p\to 1$ (resp., $p\to 0$) for the distribution of~$D$ conditional on a victory of~$A$ (resp., of~$B$).

\begin{proposition} \label{propolimas} 
Let, for the side-out scoring system, $t\mapsto M_{A}^{C}(t)={\rm E}[ e^{tD}\, |\, E^{C}, S=A ]$, $C\in\{ A,B \}$, be the moment generating function of~$D$ conditional on the event~$E^{C} \cap [S=A]$. Denote by $t\mapsto \bar{M}_{A}^{C}(t)={\rm E}[ e^{tD}\, |\, \bar{E}^{C}, S=A ]$, $C\in\{ A,B \}$,  the corresponding  moment generating function for the rally-point system. Then,
(i) as $p\to 1$, $M_\pla^\pla(t) \to e^{nt}$ and $\bar{M}_\pla^\pla(t) \to  e^{nt};$
(ii) as $p\to0$,  $M_\pla^\plb(t) \to e^{(n+1)t}$   and  $\bar{M}_\pla^\plb(t) \to  e^{nt}.$
\end{proposition}

\noindent\proof (i) By conditioning, we get $M_\pla^\pla(t) =\sum_{k=0}^{n-1} M_\pla^{n,k,\pla}(t) p_\pla^{n,k,\pla}/p_\pla^\pla$. It is easy to check that  $\lim_{p\to1} p_\pla^{n,k,\pla}/p_\pla^\pla= \delta_{k, 0}$ and that $\lim_{p\to1}M_\pla^{n,k,\pla}(t) = e^{(n+k)t}$. Hence $\lim_{p\to1}M_\pla^\pla(t) = e^{nt}$. Likewise, $\bar{M}_\pla^\pla(t) =\sum_{k=0}^{n-1} e^{(n+k)t} \bar{p}_\pla^{n,k,\pla}/\bar{p}_\pla^\pla$. Again, it is easy to check that $\lim_{p\to1}\bar{p}_\pla^{n,k,\pla}/\bar{p}_\pla^\pla= \delta_{k, 0}$. Hence, we indeed have  $\bar{M}_\pla^\pla(t) \to  e^{nt}.$ (ii) The proof is similar, and thus left to the  reader. \cqfd

\begin{corollary} \label{coro1}
(i)  As $p\to1$,  $(e_\pla^\pla,v_\pla^\pla) \to (n,0)$ and $(\bar{e}_\pla^\pla,\bar{v}_\pla^\pla) \to (n,0)$, so that, conditional on a victory of~$A$ in an $A$-game, $D\stackrel{\rm P}{\to}n$, irrespective of the scoring system; (ii) as $p\to  0$,  $(e_\pla^\plb,v_\pla^\plb ) \to (n+1,0)$ and $(\bar{e}_\pla^\plb,  \bar{v}_\pla^\plb) \to (n, 0 )$, so that, conditional on a victory of~$B$ in an $A$-game, $D\stackrel{\rm P}{\to}n+1$ (resp., $n$) for the side-out (resp., rally-point) scoring system.
\end{corollary}

As shown by Proposition~\ref{propolimas} and Corollary~\ref{coro1}, the situation is here very clear. In each of the four cases considered, only one trajectory is possible, namely that for which all rallies  in the game will be won by the winner of the game.

Next we derive the limiting conditional distribution of $D$ under events which occur with zero probability, that is, limits as~$p\to 1$ (resp., $p\to 0$) for the distribution of~$D$ conditional on a victory of~$B$ (resp., of~$A$).  Our conclusions are much more surprising. 

\begin{proposition} \label{propyvik2}
Let $m(t):=\sum_{k=0}^{n-1}e^{(n+k)t}{ n+k-1 \choose k} \big/ \sum_{k=0}^{n-1}{ n+k-1 \choose k}$. Then,  (i) as $p\to 0$, $M_\pla^\pla(t) \to e^{nt}$ and $\bar{M}_\pla^\pla(t) \to  m(t);$
(ii) as $p\to1$, $M_\pla^\plb(t) \to(e^{(n+1)t}-e^{(2n+1)t})/(n(1-e^t))$ and  $\bar{M}_\pla^\plb(t) \to  m(t)$. In particular, as $p\to1$,  the limiting   distribution of $D$  conditional on the event $E^\plb \cap [S=\pla]$ is \emph{uniform} over the set $\{n+1, \ldots, 2n\}$.
\end{proposition}

\proof  We  first prove the assertions  for the rally-point scoring system.  In this case,   $\bar{M}_\pla^\pla(t)=\sum_{k=0}^{n-1} e^{(n+k)t} \bar{p}_\pla^{n,k,\pla}/\bar{p}_\pla^\pla$. Now, from Remark \ref{amernoserve}  it is immediate that $\lim_{p\to0}p_\pla^{n,k,\pla}/p_\pla^\pla={n+k-1 \choose k} \big/ \sum_{k=0}^{n-1}{ n+k-1 \choose k}$, which proves the claim for $\bar{M}_\pla^\pla(t)$ (hence, by symmetry, also for $\bar{M}_\pla^\plb(t)$). 

Next consider  the assertions  for the side-out scoring system.  First note that, as before, $M_\pla^\pla(t) = \sum_{k=0}^{n-1}M_\pla^{n,k,\pla}(t) p_\pla^{n,k,\pla}/p_\pla^\pla$ and  $M_\pla^\plb(t) = \sum_{k=0}^{n-1}M_\pla^{k,n,\plb}(t) p_\pla^{k,n,\plb}/p_\pla^\plb$. Now fix $k \in \{0, \ldots, n-1\}$. Using Theorem \ref{probascore}, one readily shows that   $$\lim_{p\to0}{ p_\pla^{n,k,\pla}}\big/{p_\pla^\pla} = \delta_{k,0} \mbox{  and   } \lim_{p\to1}  { p_\pla^{k,n,\plb}}\big/{p_\pla^\plb} ={1}/{n}.$$ Combining these results and the definitions of the moment generating functions,  it is  then straightforward to show that   
$$\lim_{p\to0}M_\pla^{n,k,\pla}(t) = e^{(n+k)t}  \mbox{ and } \lim_{p\to1}M_\pla^{k,n,\plb}(t) = e^{(n+k+1)t}.$$ The claim follows. \cqfd

\begin{corollary} \label{coro2}
(i) As $p\to  0$,  $(e_\pla^\pla, v_\pla^\pla) \to (n,0)$ and $(\bar{e}_\pla^\pla,\bar{v}_\pla^\pla) \to (\frac{2n^2}{n+1},\frac{2n^2(n-1)}{(n+1)^2(n+2)})$; as $p\to1$,  $(e_\pla^\plb, v_\pla^\plb) \to (\frac{3n+1}{2}, \frac{(n-1)^2}{12})$ and $(\bar{e}_\pla^\plb ,\bar{v}_\pla^\plb) \to ( \frac{2n^2}{n+1},\frac{2n^2(n-1)}{(n+1)^2(n+2)})$.
\end{corollary}
\vspace{.5mm}

It is remarkable that we can again give a complete description of the ``distribution of the process" (by this, we mean that we can again list all trajectories of rallies  leading to the event considered, and give, for each such trajectory, its probability). Consider first the side-out 
scoring system. For victories of~$A$, the situation is very clear: Corollary~\ref{coro2} indeed 
\vspace{-.5mm}
yields that, 
conditional on a victory of~$A$ in an $A$-game,~$D\stackrel{\rm P}{\to}n$ as~$p\to0$, which implies that the only possible trajectory of rallies  is the one for which all rallies  in the  game are won by~$A$. Turn then to victories of~$B$. There, we obtained in the proof of Proposition~\ref{propyvik2} that all scores ($k,n$) are equally likely. It is 
\vspace{-.5mm}
actually easy to show that, conditional on~$E^{k, n, B}\cap[S=\pla]$, $D\stackrel{\rm P}{\to} n+k+1$ as~$p\to1$. This implies that there are exactly $n$ equally likely trajectories: $\pla$ first scores $k$ points, then loses his/her serve, before $\plb$ scores  $n$ (miraculous) points and wins the game ($k=0,\ldots,n-1$).  

Consider finally the rally-point system. In this case, it is sufficient to study the distribution of the scores after victories of $\pla$ Ê(when~$p\to 0$) since Êthe number of rallies  is a function of the scores only, and since Êthe conclusions will, by symmetry, Êbe Êidentical for victories of $\plb$ Ê(when~$p\to 1$). Clearly, for any fixed~$k\in\{0,1,\ldots,n-1\}$, there are exactly~${n+k-1 \choose k}$ trajectories leading to the score~$(n,k)$, and those trajectories are equally likely. Each such trajectory will then occur  with probability 
$1/\sum_{k=0}^{n-1}{ n+k-1 \choose k}$, because, as we have seen in the proof of Proposition~\ref{propyvik2}, the score~$(n,k)$ occurs with probability ${n+k-1 \choose k} \big/ \sum_{k=0}^{n-1}{ n+k-1 \choose k}$. These considerations provide the whole distribution of the process: there are $ \sum_{k=0}^{n-1}{ n+k-1 \choose k}$ equally likely possible trajectories, namely the ones we have just considered. The exact limiting distribution of~$D$ can of course trivially be computed from this.

\
\vspace{-2mm}

\bigskip
\noindent\textbf{\large Acknowledgements}
\vspace{4mm}

Davy Paindaveine, who is also  member of ECORE, the association between CORE and ECARES,  is grateful to the Fonds National de la Recherche Scientifique, Communaut\'{e} fran\c{c}aise de Belgique, for a Mandat d'Impulsion Scientifique. Yvik Swan thanks the Fonds National de la Recherche Scientifique, Communaut\'{e} fran\c{c}aise de Belgique, for support via a Mandat de Charg\'{e} de Recherches FNRS.

\newpage
\vspace{4mm} 

%




\begin{figure}[H] 
\begin{center}
\vspace{-4mm}
\includegraphics[width=10cm]{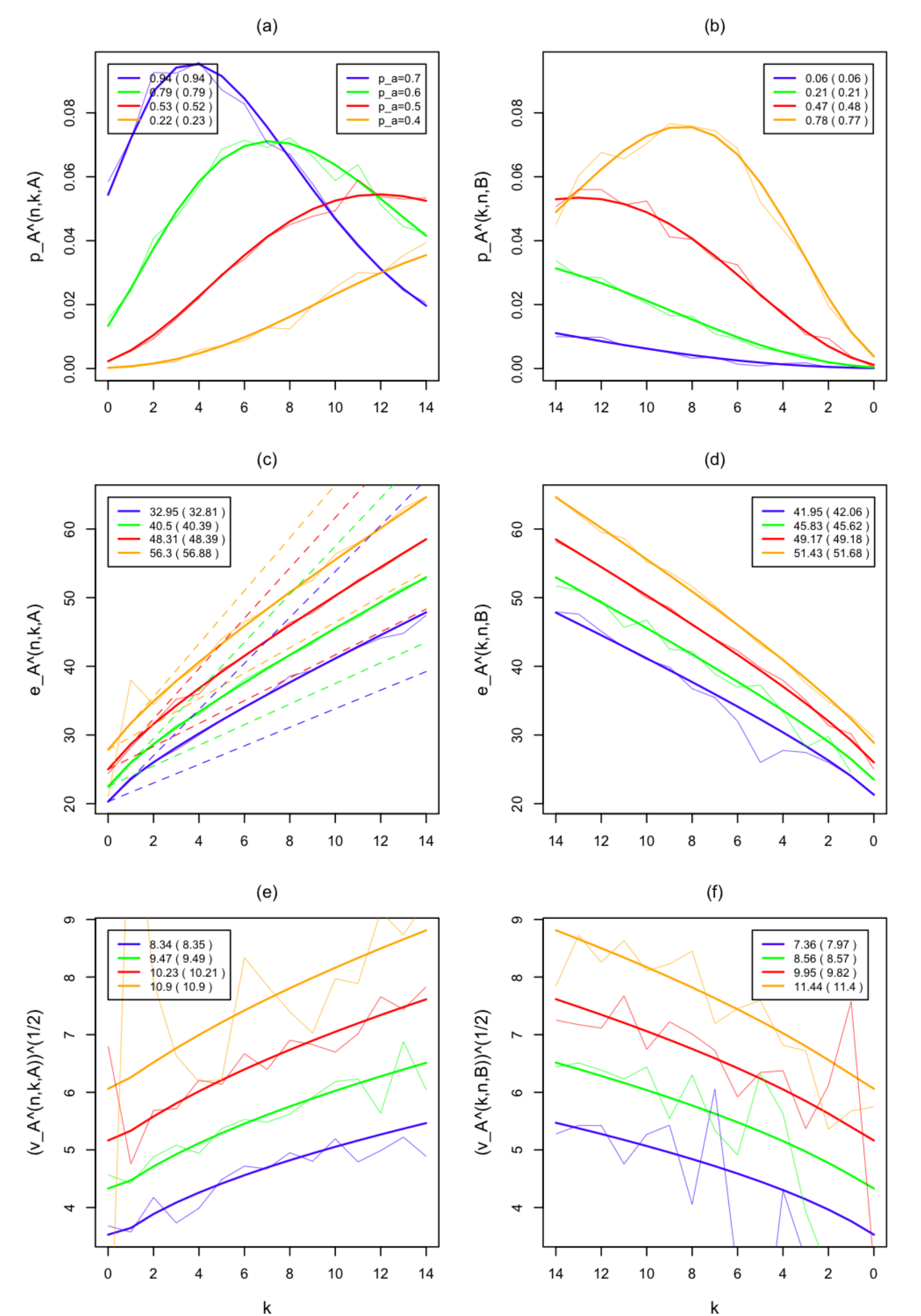}
\vspace{-6mm}
\end{center}
\caption{All subfigures refer to an $A$-game played under the  side-out scoring system with~$n=15$. Left: for~$(p_a,p_b)=(.7,.5)$, $(.6,.5)$, $(.5,.5)$, and $(.4,.5)$, (a) 	probabilities~$p_A^{n,k,A}$ that Player~$\pla$ wins the game on the score~$(n,k)$ (along with the probabilities~$p_A^{A}$ that Player~$A$ wins the game), (c) expected values~$e_A^{n,k,A}$  and (e) standard deviations~$(v_A^{n,k,A})^{1/2}$  of the numbers of rallies~$D$ conditional on the corresponding events (along with the expected values~$e_A^A$ and standard deviations~$(v_A^A)^{1/2}$ of~$D$ conditional on a victory of~$A$). Right: the corresponding values for victories of~$B$ on the score~$(k,n)$. As for the  expected values and standard deviations of~$D$ unconditional on the score or the winner, we have~$(e_A,v_A^{1/2})=	(33.5,8.6)$, $(41.6,9.5)$, $(48.7,10.1)$, and $(52.5,11.5)$, 
for~$(p_a,p_b)=(.7,.5)$, $(.6,.5)$, $(.5,.5)$, and $(.4,.5)$, respectively. Estimated probabilities, expectations, and standard deviations based on $5,000$ replications are also reported (thinner lines in plots and numbers between parentheses in legend boxes). Dashed lines in (c) correspond to Simmons' (1989) lower and upper bounds in~(\ref{simmons}).}

\label{fig1}
\end{figure}


\begin{figure}[H] 
\begin{center}
\includegraphics[width=12cm]{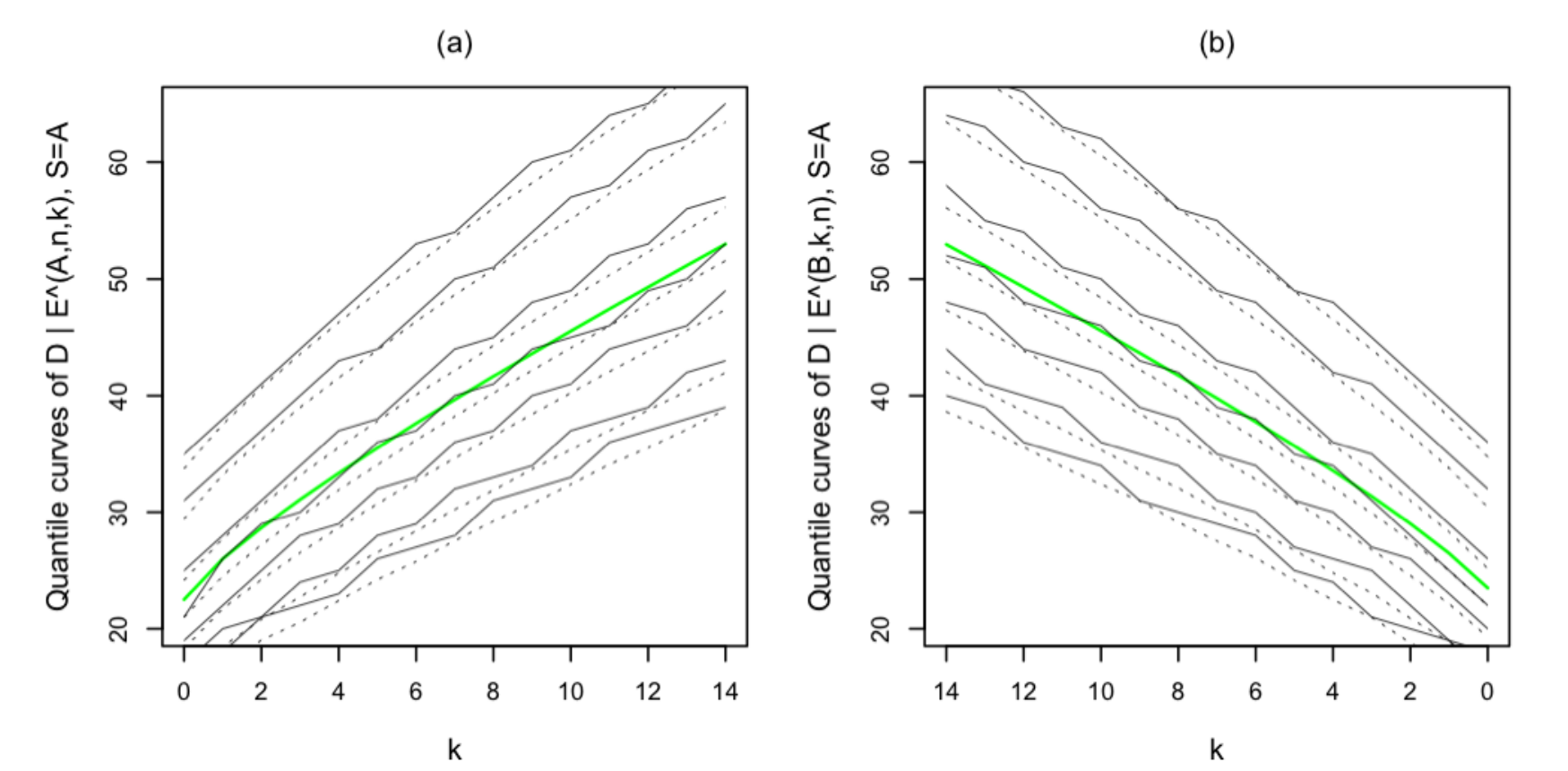}
\vspace{-5mm}
\end{center}
\caption{Both subfigures refer to an $A$-game played under the  side-out scoring system with~$n=15$ and~$(p_a,p_b)=(.6,.5)$. Subfigure~(a) (resp., Subfigure~(b)) reports, in black and as a function of~$k$, the $\alpha$-quantile of the number of rallies needed to complete the game, conditional on a victory of~$A$ on the score~$(n,k)$ (resp., conditional on a victory of~$B$ on the score~$(k,n)$), with~$\alpha=.01,.05,.25,.50,.75,.95$, and $.99$. Solid lines (resp., dotted lines) correspond to standard (resp., interpolated) quantiles; see Section~\ref{distriduration} for details. The green curves are the same as in Figure~\ref{fig1}, hence give the 	expected values of~$D$ conditional on the same events.}
\label {fig2}
\end{figure}


\begin{figure}[H] 
\begin{center}
\includegraphics[width=11cm]{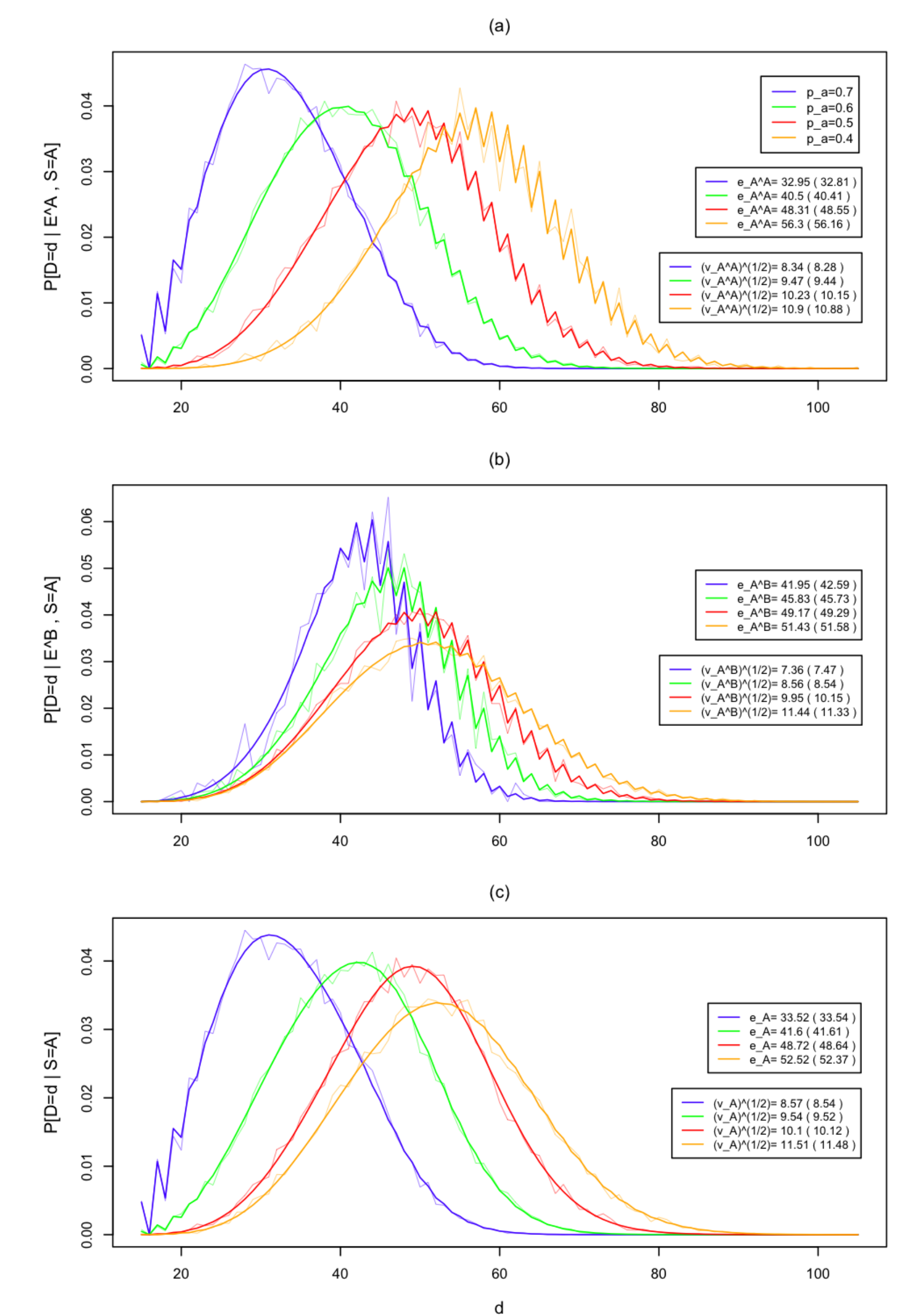}
\vspace{-5mm}
\end{center}
\caption{All subfigures refer to an $A$-game played under the  side-out scoring system with~$n=15$. For~$(p_a,p_b)=(.7,.5)$, $(.6,.5)$, $(.5,.5)$, and $(.4,.5)$, they report the probabilities that the number of rallies~$D$ needed to complete the game takes value~$d$,
(a) conditional upon a victory of Player~$A$,
(b) conditional upon a victory of Player~$B$, and
(c) unconditional. Empirical frequencies based on $20,000$ replications are also reported (thinner lines in plots and numbers between parentheses in legend boxes).}
\label {fig3}
\end{figure} 
  

\begin{figure}[H] 
\begin{center}
\includegraphics[width=13cm]{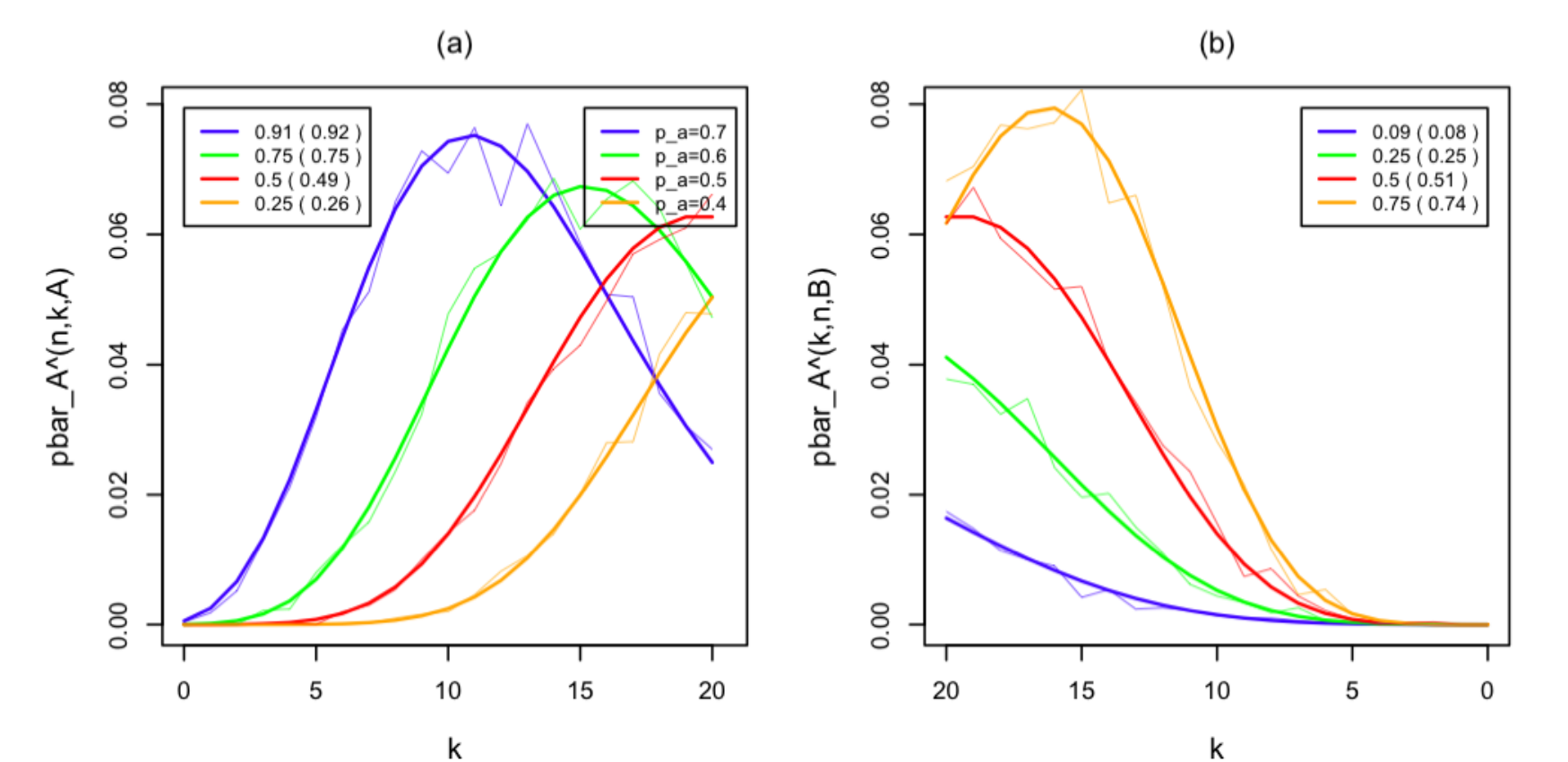}
\vspace{-5mm}
\end{center}
\caption{Both subfigures refer to an $A$-game played under the rally-point scoring system with~$n=21$. Subfigure~(a): for~$(p_a,p_b)=(.7,.5)$, $(.6,.5)$, $(.5,.5)$, and $(.4,.5)$, 
	probabilities~$\bar{p}_A^{n,k,A}$ that Player~$\pla$ wins the game on the score~$(n,k)$, along with the probabilities~$\bar{p}_A^{A}$ that Player~$A$ wins the game. Subfigures~(b): the corresponding values for victories of~$B$ on the score~$(k,n)$. Estimated probabilities based on  $5,000$ replications are also reported (thinner lines in plots and numbers between parentheses in legend boxes).}
\label{fig4}
\end{figure} 


\begin{figure}[H] 
\begin{center}
\includegraphics[width=12cm]{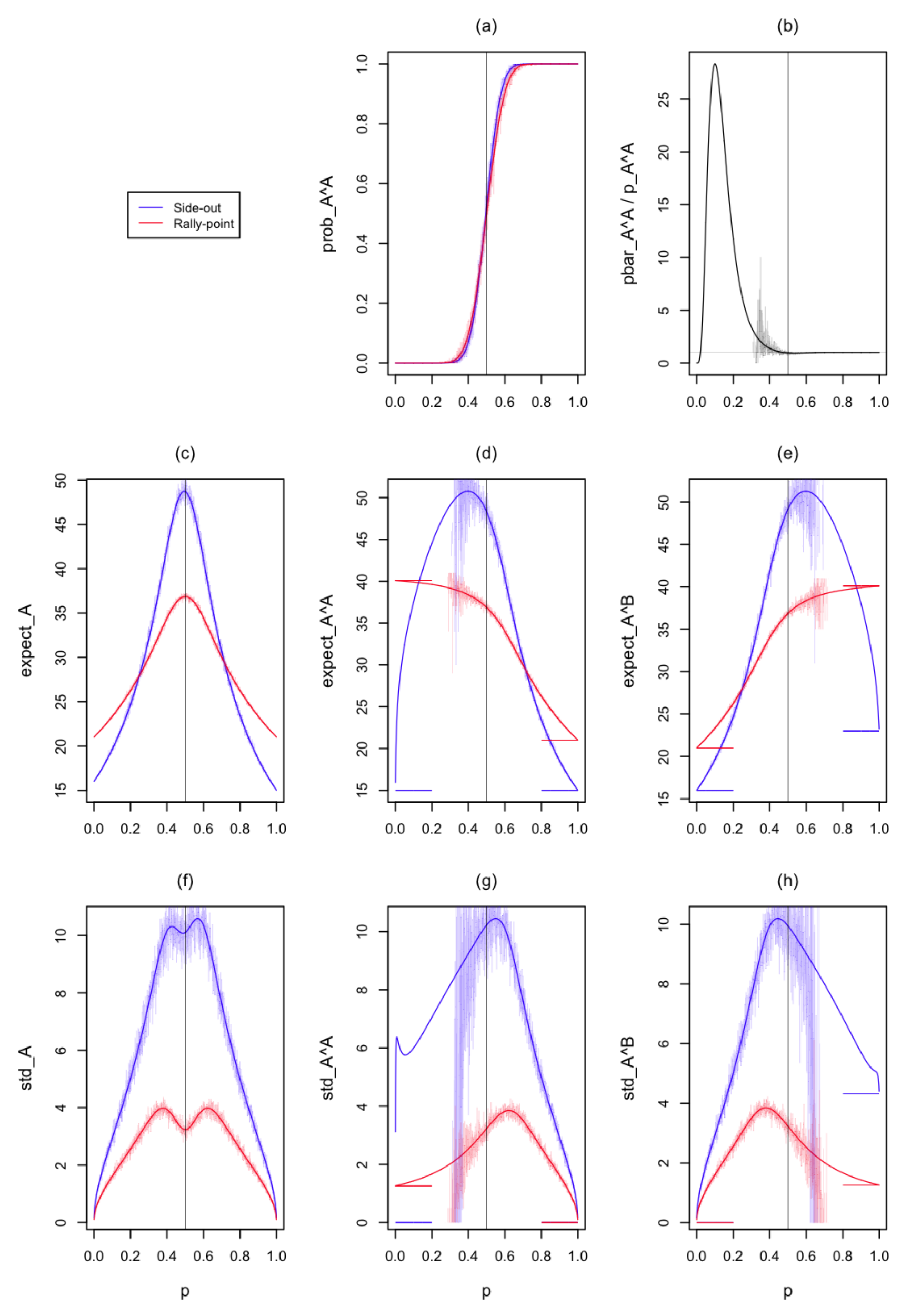}
\vspace{-8mm}
\end{center}
\caption{As a function of~$p=p_a=1-p_b$ (that is, in the no-server model), probabilities~$p_A^{15,k,A}$ (in blue) that Player~$\pla$ wins an $n=15$ side-out $A$-game on the score~$(15,k)$, along with the probabilities~$\bar{p}_A^{21,k,A}$ (in red) that Player~$\pla$ wins an $n=21$ rally-point  $A$-game on the score~$(21,k)$. Expectations (first row) and standard deviations (second row) of the number of rallies needed to complete the corresponding games, unconditional on the winner (first column),  conditional on a victory of Player~$A$ (second column),  and conditional on a victory of Player~$B$ (third column).  Estimated probabilities, expectations, and standard deviations (based on $200$ replications at each value of~$p=0,.0005,.0010,.0015,\ldots,.9995$) are also reported (thinner lines). 
}
\label{fig5}
\end{figure} 



\begin{figure}[H] 
\begin{center}
\includegraphics[width=12cm]{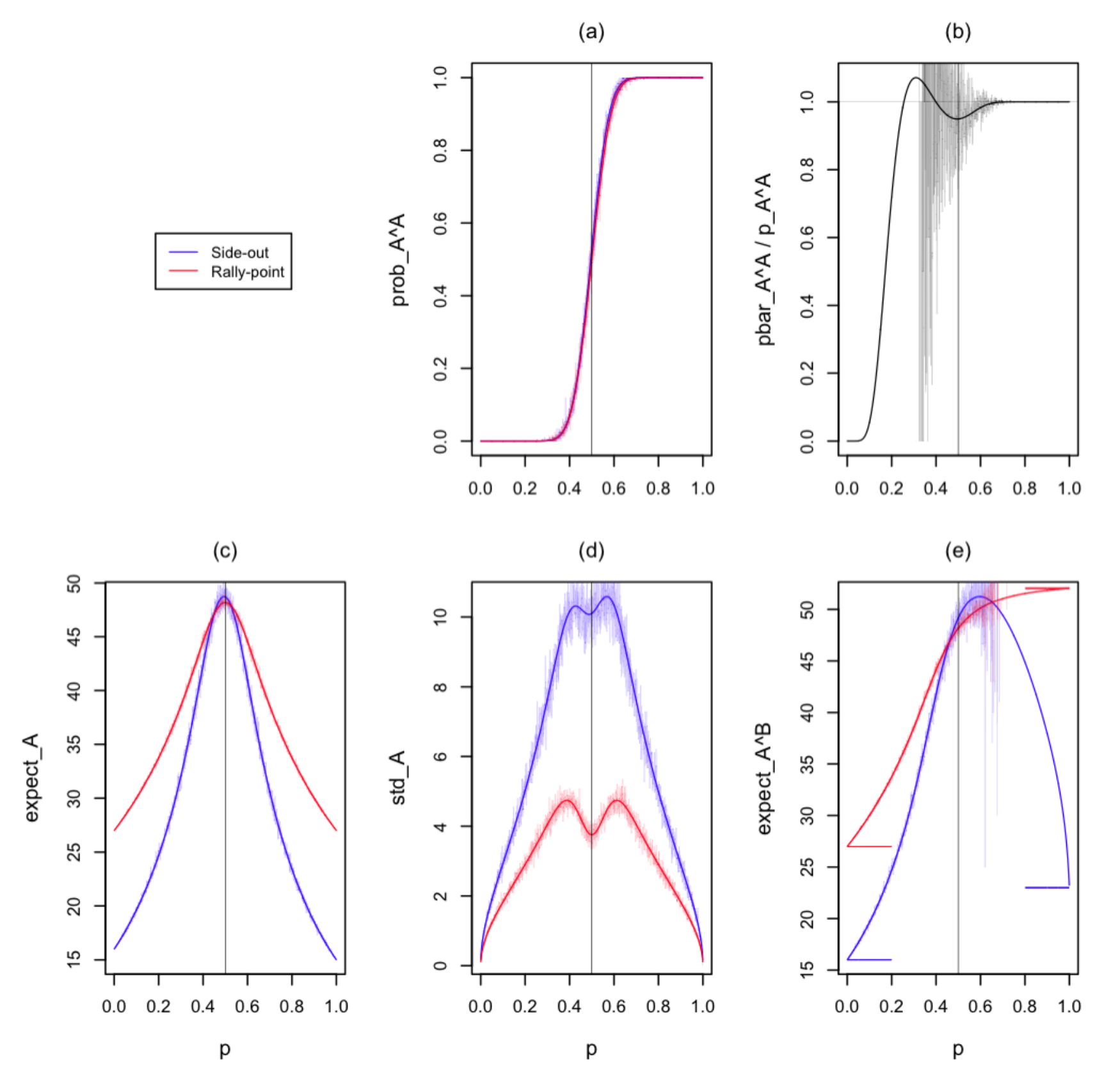}
\vspace{-4mm}
\end{center}
\caption{Subfigures (a)-(e) here report Subfigures (a)-(c), (f), and (e) from Figure~\ref{fig5} with the only difference that the rally-point scoring here is based on~$n=27$ (the side-out scoring is still based on~$n=15$).}
\label{fig6}
\end{figure}


\vspace{-10mm}
\begin{figure}[H] 
\begin{center}
\includegraphics[width=13cm]{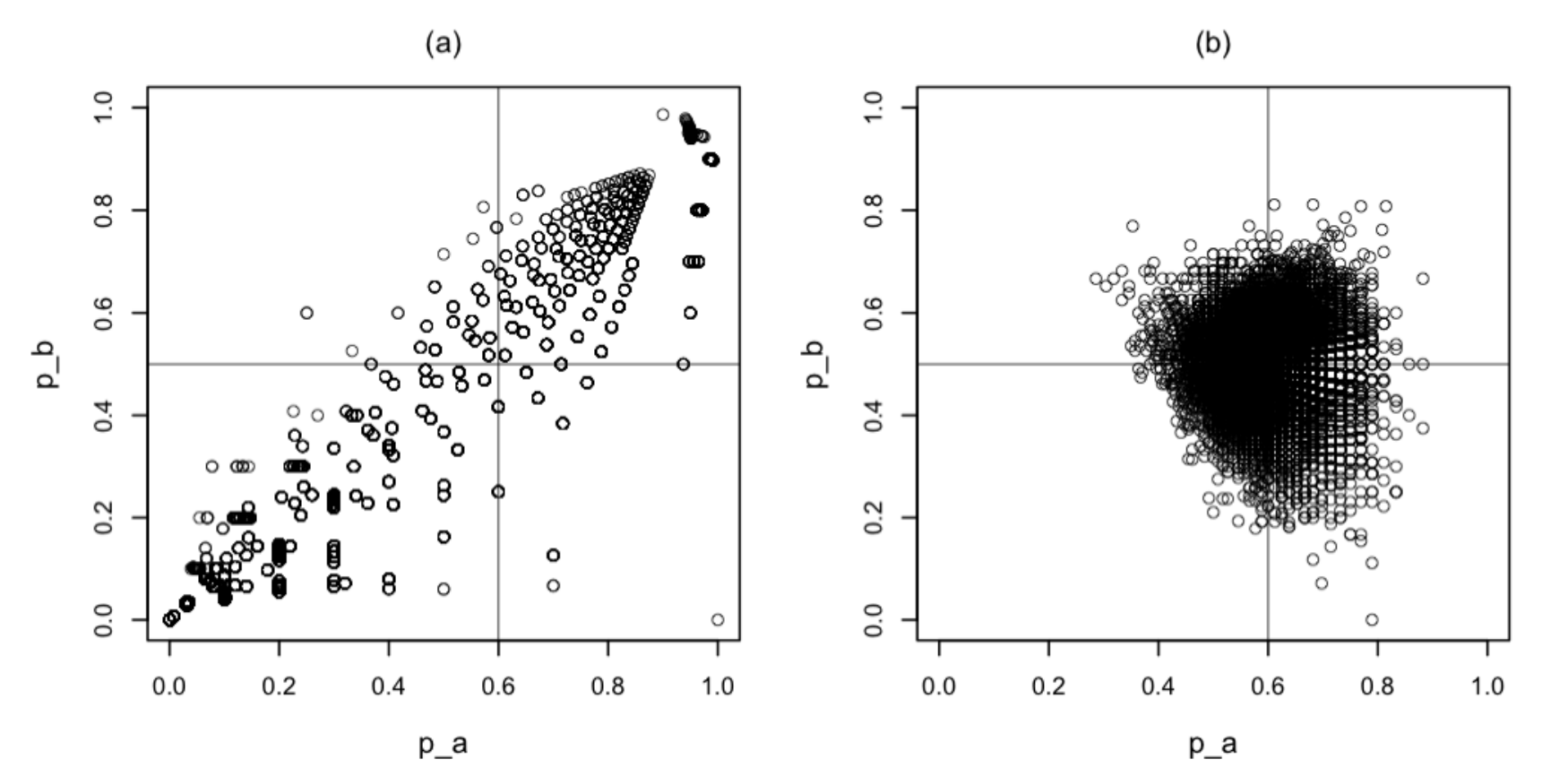}
\vspace{-5mm}
\end{center}
\caption{Subfigure~(a) (resp., (b)) is a scatter plot of the values of score-based (resp., score-and-duration-based) maximum likelihood estimators for~$(p_a,p_b)$, from  $J=1,000$ replications of $m=2$ side-out $A$-games with~$n=15$.}
\label{fig7}
\end{figure}

\end{document}